\documentclass[12pt]{amsart}
\usepackage{amssymb}
\usepackage{amsmath}
\usepackage{pictex}
\newcommand\AC{\mathcal A}
\newcommand\al{\alpha}
\newcommand\AND{\quad\mbox{and}\quad}

\newcommand\bd{\partial}
\newcommand\C{\mathbb C}
\newcommand\CC{\mathcal C}
\newcommand\cf{\curlywedge}

\newcommand\de{\delta}
\newcommand\DL{\mbox{\sl DL}}
\newcommand\dn{\mathfrak{d}}
\newcommand\dps{\displaystyle}
\newcommand\ep{\varepsilon}
\newcommand\ff{{\mathbf f}}
\newcommand\ga{\gamma}
   \newcommand\Ga{\Gamma}
\newcommand\geo[1]{\overline{#1}}

\newcommand\hor{\mathfrak{h}}
\newcommand\ka{\kappa}
\newcommand\la{\lambda}
\newcommand\lle{\preccurlyeq}
\newcommand\lra{\leftrightarrow}
\newcommand\mlt{\operatorname{\sf mult}}
   
\newcommand\om{\omega}
\newcommand\mm{\mathsf{m}}
\newcommand\N{\mathbb N}
\newcommand\nn{\mathfrak{s}}
\newcommand\R{\mathbb R}
\newcommand\scs{\scriptstyle}
\newcommand\Sig{\Sigma}
   
\newcommand\spec{\operatorname{\sf spec}}
\newcommand\spn{\operatorname{\sf span}}
\newcommand\supp{\operatorname{\sf supp}}
\newcommand\T{\mathbb T}

\newcommand\up{\mathfrak{u}}
\newcommand\wh{\widehat}
   
\newcommand\Z{\mathbb Z}

\numberwithin{equation}{section}

\newtheoremstyle{mythm}
  {9pt}
  {9pt}
  {\itshape}
  {0pt}
  {\bfseries}
  {}
  { }
  {\thmnumber{(#2)}\thmname{ #1}\thmnote{ #3}}

\newtheoremstyle{mydef}
  {9pt}
  {9pt}
  {\normalfont}
  {0pt}
  {\bfseries}
  {}
  { }
  {\thmnumber{(#2)}\thmname{ #1}\thmnote{ #3}}

\theoremstyle{mythm}
\newtheorem{thm}[equation]{Theorem.}
\newtheorem{pro}[equation]{Proposition.}
\newtheorem{lem}[equation]{Lemma.}
\newtheorem{cor}[equation]{Corollary.}

\theoremstyle{mydef}
\newtheorem{dfn}[equation]{Definition.}

\newtheorem{rmk}[equation]{Remark.}

\begin{document}
\title{\large Spectral computations on lamplighter groups and 
Diestel-Leader graphs}
\author{\bf Laurent BARTHOLDI and Wolfgang WOESS}
\address{\parbox{1.4\linewidth}{IGAT, B\^atiment BCH, 
\'Ecole Polytechnique F\'ed\'erale,\\ 
    CH-1015 Lausanne, Switzerland\\}}
\email{laurent.bartholdi@epfl.ch}
\address{\parbox{1.4\linewidth}{Institut f\"ur Mathematik C, 
Technische Universit\"at Graz,\\
Steyrergasse 30, A-8010 Graz, Austria\\}}
\email{woess@TUGraz.at}
\date{May 11, 2004}
\thanks{Supported by FWF (Austrian Science Fund) project P15577}
\begin{abstract}
The Diestel-Leader graph $\DL(q,r)$ is the horocyclic product of the homogeneous
trees with respective degrees $q+1$ and $r+1$. When $q=r$, it is the
Cayley graph of the lamplighter group (wreath product) $\Z_q \wr \Z$
with respect to a natural generating set. For the ``Simple random walk''
(SRW) operator on the latter group, {\sc Grigorchuk
and \.Zuk} and {\sc Dicks and Schick} have determined the spectrum
and the (on-diagonal) spectral measure (Plancherel measure).
Here, we show that thanks to the geometric realization, these results
can be obtained for all $\DL$-graphs by directly computing an 
$\ell^2$-complete orthonormal system of finitely supported 
eigenfunctions of the SRW. This allows computation of all matrix elements
of the spectral resolution, including the Plancherel measure. 
As one application, we determine the sharp asymptotic behaviour
of the $N$-step return probabilities of SRW. The spectral computations 
involve a natural approximating sequence of finite subgraphs, and
we study the question whether the cumulative spectral distributions of the
latter converge weakly to the Plancherel measure. To this end, we provide
a general result regarding \emph{F\o lner approximations;} in the specific case
of $\DL(q,r)$, the answer is positive only when $r=q$.
\end{abstract}
\maketitle

\markboth{\sf L. Bartholdi and W. Woess}{\sf Spectral computations}
\baselineskip 15pt

\section{Introduction}
Let $X$ be a locally finite connected graph. \emph{Simple random walk} (SRW)
on $X$ is the Markov chain on $X$ with transition probabilities
$$
p(x,y) = \begin{cases} 1/\deg(x) & \text{if}\; y \sim x\,,\\
                               0 & \text{otherwise.}
         \end{cases}
$$
Here, $\sim$ denotes neighbourhood, and $\deg(x)$ is the degree (number of
neighbours) of vertex $x \in X$. The transition operator associated with
SRW acts on real or complex functions $f$ on $X$ by
$$
Pf(x) = \sum_y p(x,y)f(y)\,.
$$
In particular, $P$ acts as a self-adjoint operator on the weighted
$\ell^2$-space $\ell^2(X,\deg)$ and has norm = spectral radius $\le 1$.
Here, we shall only consider \emph{regular} graphs,  i.e., $\deg$ is constant,
and we prefer to use the ordinary space $\ell^2(X)$, where the reference
measure is the counting measure (instead of $m(x) = \deg(x)$).

Associated with $P$ there is the \emph{resolution of the identity.}
This is an operator-valued measure $E$ defined on $\spec(P)$ such that
for all $N \ge 0$,
$$
P^N = \int_{\spec(P)}  \la^n E(d\la)
$$
It is characterized by the matrix elements
\begin{equation}\label{measures}
\mu_{x,y}(d\la) = \langle \de_x,E(d\la)\de_y\rangle\,,
\end{equation}
which in turn are characterized by their moments, which are the $N$-step
transition probabilities (matrix elements of $P^N$),
\begin{equation}\label{moments}
p^{(N)}(x,y) = \int_{\spec(P)}  \la^N\, \mu_{x,y}(d\la) \quad \forall\ N \ge 0
\,.
\end{equation}
Each $\mu_{x,x}$ is a probability measure, while the off-diagonal
$\mu_{x,y}$ are signed measures with total mass $0$.
When $X$ is a \emph{transitive} graph (i.e., its group of isometries
acts transitively on the vertex set), then all $\mu_{x,x}$ coincide,
and we shall just write $\mu$ for this measure, whose support is the whole
spectrum. This holds, in particular, for \emph{Cayley graphs} of 
finitely generated groups.\footnote{If $\Ga$ is a group and 
$S=S^{-1} \subset \Ga$ a finite set of generators, then the Cayley graph
$X(\Ga,S)$ has vertex set $\Ga$, and $x \sim y$ if $x^{-1}y \in S$.} 
In the spirit of Harmonic Analysis, we call $\mu$
the \emph{Plancherel measure;} more recently, it has also been
called the \emph{Kesten spectral measure} by some authors.

Basic references for the general theory of spectra of \emph{infinite} 
graphs and groups are the papers of {\sc \mbox{Mohar} and Woess} \cite{MoWo}, 
{\sc de la Harpe, Robertson and Valette} \cite{HaRoVa} and 
{\sc \mbox{Grigorchuk} and
\.Zuk} \cite{GrZu1}. In \cite{MoWo} and \cite{HaRoVa}, one can also find
many examples of specific graphs and groups where spectra and (less
frequently) spectral measures are computed: in basic cases 
$\spec(P)$ is an interval,
and the Plancherel measure has a continuous density with respect to Lebesgue measure.
This occurs for integer lattices -- a classical result from Fourier analysis, 
see e.g. {\sc P\'olya} \cite{Po} -- and for free groups, resp. homogeneous
trees -- see {\sc Kesten} \cite {Ke} and {\sc Cartier} \cite{Ca}.
For other tree-like cases (infinite distance-regular graphs), the situation
is almost the same, with a possible additional isolated eigenvalue, 
as was shown by {\sc Kuhn and Soardi} \cite{KuSo}; see also 
{\sc Faraut and Picardello} \cite{FaPi}. For $\tilde A_d$-buildings with
arbitrary $d$, the situation is similar to that of homogeneous trees,
see {\sc Cartwright and M\l otkowski} \cite{CarMl} and {\sc Cartwright} 
\cite{Car}.

The situation is different on typical fractal graphs such as
the one associated with the Sierpi\'nski gasket, where the spectrum
is \emph{pure point}, i.e., the closure of the set of eigenvalues of $P$,
see {\sc Malozemov and Teplyaev} \cite{MaTe}, {\sc Teplyaev}
\cite{Te}, {\sc Sabot} \cite{Sa} and, for a 
generalization, {\sc Kr\"on} \cite{Kr}. These graphs are regular, but far
from being transitive. However, it turned out in recent research
that a similar situation may also occur in certain classes of \emph{fractal
groups} related with the construction of {\sc Grigorchuk} and 
{\sc Gupta and Sidki} of finitely generated infinite torsion groups with
intermediate growth. For a comprehensive survey of these groups and their 
properties, see {\sc Bartholdi, Grigorchuk and Nekrashevych} \cite{BaGrNe}
and the references given there,
and for the specific computation of a pure point spectrum on such a group,
see {\sc Bartholdi and Grigorchuk} \cite{BaGr}. 

Coming finally to the types of structures considered in the present
paper, a pure point spectrum with the associated Plancherel (Kesten)
measure was recently detected for a different class of groups, namely
the \emph{lamplighter groups} (wreath products) $F \wr \Z$,
  where $F$ is a finite group; see {\sc Grigorchuk and \.Zuk}
\cite{GrZu2} for $\Z_2 \wr \Z$, and {\sc Dicks and Schick}
\cite{DiSc} for the general case\footnote{We denote by $\Z_q$ 
the cyclic group $\Z/q\Z$ of order $q$}. While \cite{GrZu2} uses
approximation of the considered Cayley graph by an increasing sequence
of finite graphs, \cite{DiSc} applies von-Neumann-algebraic methods.

Note that the present study of random walks on $F\wr\Z$ depends only 
on the cardinality $|F|$, whence it suffices to consider $F=\Z_q$.
Now, the Cayley graphs of $\Z_q \wr \Z$ considered in \cite{GrZu2} and
\cite{DiSc} turn out to have geometric realizations as specific 
examples $\DL(q,q)$ in the family of \emph{Diestel-Leader graphs}
introduced by {\sc Diestel and Leader} \cite{DiLe},
i.e., the horocyclic products $DL(q,r)$ of two homogeneous trees
with degree $q+1$ and $r+1$, respectively; see {\sc Woess}
\cite{Wo} and \S \ref{geometry} below. 

In the present paper, we exploit this geometric model to provide in \S
\ref{tetrahedra} a completely explicit and elementary construction of
an $\ell^2$-complete orthonormal system of finitely supported
eigenfunctions of the SRW-operator on $DL(q,r)$.  This comprises the
lamplighter groups, but holds more generally for all Diestel-Leader
graphs. When $q \ne r$, the graph $DL(q,r)$ is transitive, but not a
Cayley graph, whence the group-specific methods of \cite{DiSc} do not
apply here (they rely on identifying eigenfunctions as
projections in the von Neumann algebra of a group acting on the
graph; here that group is non-discrete and non-unimodular), 
nor can one use an approximating sequence of Schreier graphs as in \cite{GrZu2}.

We recover the spectral radius of SRW on $DL(q,r)$, computed
previously by {\sc Saloff-Coste and Woess} \cite{SCWo},
\begin{equation}\label{specrad}
\rho(P) = \frac{2\sqrt{qr}}{q+r}\,.
\end{equation}
The spectrum $\spec(P) = [-\rho(P)\,,\,\rho(P)]$ is the closure of the
set of all eigenvalues
$$
\left\{\la_{m,n} = \rho(P) \cos \frac{m}{n}\pi : n \ge 2\,,\;1 \le m \le n-1
\right\},
$$
see Theorem \ref{spec}. This theorem can be used to find an expression for each
of the spectral measures, see \S \ref{meas}. In particular, the Plancherel
measure can be computed explicitly (Corollary \ref{muoo}). 

We then use the latter in \S \ref{asymptotics} 
to determine the sharp asymptotic behaviour
of the $2N$-step return probabilities $p^{(2N)}(x,x)$, as $N \to \infty$,
see Theorem \ref{return} (note that $p^{(2N+1)}(x,x)=0$ since
$DL(q,r)$ is bipartite). 
For the lamplighter groups, i.e., on $DL(q,q)$,
these asymptotics have been determined for almost the same random walk
by {\sc Revelle} \cite{Re}; one has
$$
p^{(2N)}(x,x) \sim 
\bar A_1 \, \exp\bigl(-B_1\,N^{1/3}\bigr)\, N^{1/6}\qquad (r=q)\,.
$$ 
It is interesting to note that in the case $r \ne q$, the last (polynomial) 
term changes by a factor of $N\,$:
$$
p^{(2N)}(x,x) \sim 
A_1\,\rho(P)^{2N} \,
\exp\bigl(-B_1\,N^{1/3}\bigr)\, N^{-5/6}
\qquad (r \ne q)\,.
$$
Furthermore, the constants $\bar A_1, A_1, B_1 > 0$ are determined
explicitly as functions of $q$ and~$r$.

Next, in \S \ref{cumulative}, we discuss for general vertex transitive 
graphs under which conditions the cumulative spectral measures of
an approximating sequence of finite graphs converge (resp. do not converge)
weakly to the Plancherel measure. A positive criterion is given in
terms of \emph{F\o lner sequences} of approximating subgraphs, and more
generally, \emph{F\o lner approximations}  -- see Theorem \ref{expanding} and
Remark \ref{adjacency}. In our setting, this applies to $\DL(q,q)$
with the natural subgraphs \emph{(tetrahedra)} used in \S \ref{tetrahedra} 
for computing the spectrum, since they constitute a F\o lner sequence.
However, when $r \ne q$, this is not true, and the cumulative
spectral measures associated with tetrahedra do not converge to
the Plancherel measure -- see Proposition \ref{tildemun}, which is 
preceded by lengthy computations of the spectra of tetrahedra.
This should be seen in the light of \emph{amenability}. A graph is called 
amenable if it has a F\o lner sequence of subgraphs, which is equivalent with
$\rho(P) = 1$ for SRW. By \cite{SCWo}, $\DL(q,r)$ is amenable if and only
if $r=q$; see \eqref{specrad} above.
The discussion of \S \ref{cumulative} should be compared with the
results of {\sc Serre} \cite{Ser}, who studies (among other) the
question under which conditions the cumulative spectral measures
of an arbitrary sequence of \emph{regular} graphs have a weak limit.
(Our tetrahedra are not regular at their boundaries.)

At the end, in \S \ref{final}, we add several observations, including
random walks with drift, the corresponding spectra, their  
return probabilities, and also their projections on the two subtrees and
on $\Z$.

\section{Lamplighter groups and Diestel-Leader graphs}\label{geometry}

This section is a short version of \S 2 in \cite{Wo}. We explain the structure
of the DL-graphs and their relation with the groups $\Z_q \wr \Z$.

Let $\T = \T_q$ be the homogeneous tree with degree $q+1$, $q \ge 2$.
A \emph{geodesic path}, respectively \emph{geodesic ray}, respectively \emph{infinite
geodesic} in $\T$ is a finite, 
respectively one-sided infinite, respectively doubly infinite sequence $(x_n)$ of vertices
of $\T$ such that $d(x_i,x_j) = |i-j|$ for all $i, j$, 
where $d(\cdot,\cdot)$ denotes the graph distance. 

Two rays are \emph{equivalent} if the symmetric difference of their
supports is finite.
An \emph{end} of $\T$ is an equivalence class of rays. The space of 
ends is denoted $\bd \T$, and we write $\wh \T = \T \cup \bd \T$. 
For all $w, z \in \wh \T$ there is a unique geodesic $\geo{w\,z}$ 
that connects the two. In particular, if $x \in \T$ and $\xi \in \bd \T$ then 
$\geo{x\,\xi}$ is the ray that starts at $x$ and represents $\xi$.

For $x,y \in \T$, $x \ne y$, we define the cone
$\wh \T(x,y) = \{ w \in \wh \T : y \in \geo{x\,w} \}$.
The collection of all cones is the basis of  a topology which 
makes $\wh \T$ a compact, totally disconnected Hausdorff space
with $\T$ as a dense, discrete subset. 

We fix a root  $o \in \T$.
If $w, z \in \wh \T$, then their \emph{confluent} $c=w \wedge z$ with
respect to the root vertex $o$ is defined by
$\geo{o\,w} \cap \geo{o\,z} = \geo{o\,c}$. 
Similarly, we choose and fix a reference end $\om \in \bd \T$. For 
$z, v \in \wh \T \setminus \{ \om \}$, their confluent $b = v \cf z$ 
with respect to $\om$ is defined by
$\geo{v\,\om} \cap \geo{z\,\om} = \geo{b\,\om}$. We write
$$
z \lle v \quad\mbox{if}\quad z \cf v = z\,.
$$ 
For $x \in \T$, we describe its
relative position with respect to $o$ by the two numbers
$$
\up(x) = d(o,x \cf o) \AND \dn(x) = d(x,x \cf o)\,.
$$
In Figure 1, $\up(x)$ and $\dn(x)$ correspond to the numbers of steps one 
has to take upwards (in direction of $\om$), respectively downwards, on the geodesic
path from $o$ to $x$. Thus, $d(x,o) = \up(x)+\dn(x)$.

The \emph{Busemann function} $\hor: \T \to \Z$ and the \emph{horocycles} $H_k$
with respect to $\om$ are defined as
$$
\hor(x) = \dn(x) - \up(x) \AND H_k = \{ x \in \T : \hor(x) = k \}\,.
$$
Every horocycle is infinite. Every vertex $x$ in $H_k$ has one neighbour 
$x^-$ (its predecessor) in $H_{k-1}$ and $q$ neighbours (its successors)
in $H_{k+1}$. We set $\bd^* \T = \bd \T \setminus \{\om\}$. 

We label each edge of $\T$ by an element of $\Z_q$ such that for each vertex,
the ``downward'' edges to its $q$ successors carry labels $0, \dots, q-1$
from left to right (say), see Figure 1.  
Thus, for each $x \in  \T$, the sequence $\bigl( \sigma(n) \bigr)_{n \le 0}$
of labels on the geodesic  $\geo{x\,\om}$ has 
finite support $\{ n : \sigma(n) \ne 0 \}$. 
We write $\Sigma_q$ for the set of all those sequences.
On every horocycle, there is exactly one vertex corresponding to each 
$\sigma \in \Sigma_q$.
Thus, $\T_q$ is in one-to-one correspondence with the set 
$\Sigma_q \times \Z$, and the $k$-th horocycle is 
$H_k = \Sigma_q \times \{k\}$.

\vspace{-.4cm}

$$
\beginpicture 

\setcoordinatesystem units <.7mm,1.04mm>

\setplotarea x from -10 to 104, y from 4 to 84

\arrow <6pt> [.2,.67] from 2 2 to 80 80

\plot 32 32 62 2 /

 \plot 16 16 30 2 /

 \plot 48 16 34 2 /

 \plot 8 8 14 2 /

 \plot 24 8 18 2 /

 \plot 40 8 46 2 /

 \plot 56 8 50 2 /

 \plot 4 4 6 2 /

 \plot 12 4 10 2 /

 \plot 20 4 22 2 /

 \plot 28 4 26 2 /

 \plot 36 4 38 2 /

 \plot 44 4 42 2 /

 \plot 52 4 54 2 /

 \plot 60 4 58 2 /



 \plot 99 29 64 64 /

 \plot 66 2 96 32 /

 \plot 70 2 68 4 /

 \plot 74 2 76 4 /

 \plot 78 2 72 8 /

 \plot 82 2 88 8 /

 \plot 86 2 84 4 /

 \plot 90 2 92 4 /

 \plot 94 2 80 16 /


\setdots <3pt>
\putrule from -4.8 4 to 102 4
\putrule from -4.5 8 to 102 8
\putrule from -2 16 to 102 16
\putrule from -1.7 32 to 102 32
\putrule from -1.7 64 to 102 64
\setdashes <2pt>
\putrule from -1.7 -8 to 102 -8

\put {$\vdots$} at 32 -2
\put {$\vdots$} at 64 -2

\put {$\dots$} [l] at 103 6
\put {$\dots$} [l] at 103 48

\put {$H_{-3}$} [l] at -13 64
\put {$H_{-2}$} [l] at -13 32
\put {$H_{-1}$} [l] at -13 16
\put {$H_0$} [l] at -13 8
\put {$H_1$} [l] at -13 4
\put {$\bd^* \T$} [l] at -13 -8
\put {$\vdots$} at -10 -3
\put {$\vdots$} [B] at -10 70
\put {$\circ$} at 8 8
\put {$\omega$} at 82 82

\put {\scriptsize $0$} at 3.6 6.2
\put {\scriptsize $1$} at 12.2 6.2
\put {\scriptsize $0$} at 19.8 6.2
\put {\scriptsize $1$} at 28.4 6.2
\put {\scriptsize $0$} at 36   6.2
\put {\scriptsize $1$} at 44.2 6.2
\put {\scriptsize $0$} at 51.8 6.2
\put {\scriptsize $1$} at 60   6.2
\put {\scriptsize $0$} at 67.6 6.2
\put {\scriptsize $1$} at 76   6.2
\put {\scriptsize $0$} at 83.8 6.2
\put {\scriptsize $1$} at 92.2 6.2

\put {\scriptsize $0$} at 9 12
\put {\scriptsize $1$} at 22.5 12
\put {\scriptsize $0$} at 41 12
\put {\scriptsize $1$} at 54.5 12
\put {\scriptsize $0$} at 73 12
\put {\scriptsize $1$} at 86.5 12

\put {\scriptsize $0$} at 21 24
\put {\scriptsize $1$} at 43 24
\put {\scriptsize $0$} at 85 24

\put {\scriptsize $0$} at 45 48
\put {\scriptsize $1$} at 83 48

\put {\scriptsize $0$} at 72 75

\endpicture
$$

\vspace{.1cm}

\begin{center}
\centerline\emph{Figure 1}
\end{center}

\vspace{.1cm}

Now consider two trees $\T^1=\T_q$ and $\T^2=\T_r$ with roots $o_1$ 
and $o_2$ and reference ends $\om_1$ and $\om_2$, respectively.

\begin{dfn}\label{DLdef}
The Diestel-Leader graph $\DL(q,r)$ is
$$
\DL(q,r) = \{ x_1x_2 \in \T_q \times \T_r : \hor(x_1)+\hor(x_2) = 0 \}\,,
$$
and neighbourhood is given by
$\;
x_1x_2 \sim y_1y_2 \iff x_1 \sim y_1 \AND x_2 \sim y_2\,.
$
\end{dfn}

To visualize $\DL(q,r)$, draw $\T_q$ in horocyclic layers with
$\om_1$ at the top and $\bd^*\T_q$ at the bottom, and right to it $\T_r$
in the same way, but upside down, with the respective horocycles $H_k(\T_q)$
and $H_{-k}(\T_r)$ aligned. Connect the two origins $o_1$, $o_2$ by
an elastic spring. It is allowed to move along each of the two trees, 
may expand infinitely, but must always remain in horizontal position. 
The vertex set of $\DL_{q,r}$ consists of all admissible positions of 
the spring. From a position
$x_1x_2$ with $\hor(x_1) + \hor(x_2) =0$ the spring may
move downwards to one of the $r$ successors of $x_2$ in $\T_r$, and at the same 
time to the predecessor of $x_1$ in $\T_q$, or it may move upwards 
in the analogous way. Such a move corresponds
to going to a neighbour of $x_1x_2$. We see that $\DL(q,r)$ is regular
with degree $q+r$.
As the reference point in $\DL(q,r)$, we choose $o=o_1o_2$.
Figure 2 illustrates $\DL(2,2)$.

The relative position of $x=x_1x_2 \in \DL(q,r)$ with respect to $o$
is given by the four numbers $\up(x_1), \dn(x_1),  \up(x_2), \dn(x_2)$,
which satisfy the relation
\begin{equation}\label{updown}
\up(x_1) +  \up(x_2)=\dn(x_1) + \dn(x_2)\,.
\end{equation}

$$
\beginpicture
\setcoordinatesystem units <3mm,3.5mm> 

\setplotarea x from -4 to 30, y from -3.8 to 6.4
\arrow <5pt> [.2,.67] from 4 4 to 1 7
\put{$\omega_1$} [rb] at 1.2 7.2

\put{$o_1$} [lb] at  8.15 0.2

\plot -4 -4       4 4         /         
\plot 4 4         12 -4          /      
\plot -2 -2       -2.95 -4 /            
\plot -.5 -2      -1.9 -4     /         
\plot -.5 -2      -.85 -4   /           
\plot 1 -2        .2 -4      /          
\plot 1 -2        1.25  -4     /        
\plot 2.5 -2      2.3  -4     /         
\plot 2.5 -2      3.35 -4      /        
\plot 5.5 -2      4.65  -4    /         
\plot 5.5  -2     5.7  -4      /        
\plot 7  -2       6.75   -4   /         
\plot 7 -2        7.8  -4      /        
\plot 8.5  -2     8.85  -4    /         
\plot 8.5  -2     9.9  -4      /        
\plot 10  -2      10.95  -4   /         
\plot 0  0        -.5  -2      /        
\plot 2 0         1 -2     /            
\plot 2 0         2.5   -2     /        
\plot 6 0         5.5    -2    /        
\plot 6 0         7 -2         /        
\plot 8 0         8.5 -2       /        
\plot 2 2         2 0         /         
\plot 6 2         6 0         /         

\arrow <5pt> [.2,.67] from 22 -4 to 25 -7
\put{$\omega_2$} [lt] at 25.2 -7.2

\put{$o_2$} [rt] at  17.95 -.2

\plot 14  4       22 -4       /         
\plot 22 -4         30 4         /      
\plot 16 2       15.05 4 /              
\plot 17.5 2     16.1  4     /          
\plot 17.5 2     17.15  4   /           
\plot 19  2       18.2  4      /        
\plot 19 2         19.25   4     /      
\plot 20.5 2       20.3   4     /       
\plot 20.5  2      21.35 4      /       
\plot 23.5 2       22.65  4    /        
\plot 23.5  2      23.7   4      /      
\plot 25   2       24.75   4   /        
\plot 25  2        25.8  4      /       
\plot 26.5   2     26.85  4    /        
\plot 26.5   2     27.9   4      /      
\plot 28   2      28.95   4   /         
\plot 18  0        17.5  2      /       
\plot 20 0         19  2     /          
\plot 20 0         20.5   2     /       
\plot 24 0         23.5    2    /       
\plot 24 0         25 2         /       
\plot 26 0         26.5 2       /       
\plot 20 -2        20 0         /       
\plot 24 -2        24 0         /       
\put {$\circ$} at 8 0
\put {$\circ$} at 18 0
\plot 8.25 0  12.1 0 /
\plot 13.9 0 17.78 0 /
\plot    12.1   0    12.25 .4    12.25 -.4   12.55 .4   12.55 -.4
        12.85 .4   12.85 -.4  13.15 .4   13.15 -.4  13.45 .4
        13.45 -.4  13.75 .4   13.75 -.4  13.9 0     13.9  0 /

\setdashes <2pt>
\putrule from -4.5 -7  to  12.5 -7
\putrule from  13.5 7  to  30.5 7

\put {$\bd^*\T_q$} [r] at -5 -7
\put {$\bd^*\T_r$} [l] at 31 7

\put {$\vdots$} at 4 -5.2
\put {$\vdots$} at 22 5.5

\endpicture
$$
\vspace{.4cm}

\begin{center}
\emph{Figure 2}
\end{center}

\vspace{.4cm}

The \emph{lamplighter group} $\Z_q \wr \Z$ is defined as follows:
Consider the group of all finitely supported \emph{configurations} 
$$
\CC = \{ \eta: \Z \to \Z_q \,,\quad |\supp(\eta)| < \infty  \}
$$
with pointwise addition modulo $q$. Then the group $\Z$ acts on $\CC$
by translations $k \mapsto T_k:\CC \to \CC$ with $T_k\eta(m) = \eta(m-k)$.
The resulting semidirect product $\Z \rightthreetimes \CC$ is
$$
\Z_q \wr \Z = \{ (\eta,k) : \eta \in \CC\,,\;k \in \Z\} 
\quad \mbox{with group operation}\quad
(\eta,k)(\eta',k') = (\eta + T_k\eta', k+k')
$$
We identify each $(\eta,k) \in \Z_q \wr \Z$ with the vertex 
$x_1x_2 \in \DL(q,q)$, where according to the identification 
$\T_q \lra \Sigma_q \times \Z$,
the vertices $x_i$  are given by 
\begin{equation}\label{identif}
\begin{gathered}
x_1 = (\eta_k^-,k) \AND x_2=(\eta_k^+,-k)\,,\quad\mbox{where}\\
\eta_k^- = \eta|_{(-\infty\,,\,k]} \AND 
\eta_k^+ = \eta|_{[k+1\,,\,\infty)}\,,
\end{gathered}
\end{equation}
both written as sequences over the non-positive integers.

This is clearly a one-to-one correspondence between $\Z_q \wr \Z$ and 
$\DL(q,q)$, and it is also straightforward that this group acts 
transitively and fixed-point-freely on the graph.
The action of $m \in \Z$ is given by $x_1x_2 = (\sigma_1,k)(\sigma_2,-k) 
\mapsto y_1y_2 = (\sigma_1,k+m)(\sigma_2,-k+m)$,
and the action of the group of configurations is pointwise addition modulo $q$. 
Write $\de_k^{\ell}$ for the configuration in $\CC$ with value 
$\ell$ at $k$ and $0$ elsewhere.
Then $\DL(q,q)$ is the (right) Cayley graph of $\Z_q \wr \Z$ with respect to
the symmetric set of generators
$$
\{ (\de_1^{\ell},1)\,,\; (\de_0^{\ell},-1) : \ell \in \Z_q \}\,,
$$
i.e., an edge corresponds to multiplying with a generator on the right.
This is precisely the set of generators considered in \cite{GrZu2}
and \cite{DiSc} when computing the spectrum of the associated 
SRW-operator.

\medskip

\section{Tetrahedra and horizontal functions}\label{tetrahedra}

In the sequel, we shall often write $\DL$ for $\DL(q,r)$.
We say that a function $f: \DL \to \C$ is \emph{horizontal} if
it is finitely supported and
\begin{equation}\label{horizontal}
\sum_{y_2 \in \T^2 : \hor(y_2) = -\hor(x_1)} f(x_1y_2) 
= \sum_{y_1 \in \T^1 : \hor(y_1) = -\hor(x_2)} f(y_1x_2) = 0 
\qquad \forall\  x_1 \in \T^1\,,\;x_2 \in \T^2\,.
\end{equation}

\begin{lem}\label{dense}
The linear space of horizontal functions is dense in $\ell^2(\DL)$.
\end{lem}

\begin{proof} It is sufficient to show that every point mass can be approximated
in the $\ell^2$-norm by horizontal functions. Furthermore, by
vertex-transitivity, it is sufficient to show this for $\de_o$ where $o=o_1o_2$.
Let $b_1 = b_1^n \in \T^1$ be a vertex on $H_{-n}^1$
(horocycle in $\T^1$) for which $\up(b_1)=n+1$ and $\dn(b_1)=1$.
Define a function $f_1=f_1^n$ on $\T^1$ by
$$
f_1(x_1) = \begin{cases} 1 & \text{if}\; x_1 = o_1\,,\\
                        -1/q^n &\text{if}\; b_1 \lle x_1 \in H_0^1\,,\\
                         0 & \text{otherwise.}
           \end{cases}
$$
In the same way, but replacing $q$ with $r$, 
we define a function $f_2=f_2^n$ on $\T^2$.
Then the function $f=f_n$, given by $f(x_1x_2) = f_1(x_1)f_2(x_2)$,
is horizontal, and
$$
\hspace*{2.7cm}\| f_n - \de_o\|^2 = \frac{1}{q^nr^n} + \frac{1}{q^n} + \frac{1}{r^n} \to 0
\qquad\mbox{as}\qquad n \to \infty  \hspace{2.7cm}\qedhere
$$
\end{proof}

\begin{dfn}\label{tetra}
Let $a_1 \in \T^1$ and $a_2 \in \T^2$ be two vertices with
$-\hor(a_2) = \hor(a_1) + n$, where $n \ge 0$. Then the
(induced) subgraph of $\DL$ given by
$$
S = S(a_1,a_2) = \{ x_1x_2 \in \DL : a_1 \lle x_1\,,\;a_2 \lle x_2 \}
$$
is called a \emph{tetrahedron} in $\DL$ with \emph{height} $n(S)=n$.
\end{dfn}
$$
\beginpicture
\setcoordinatesystem units <3mm,3.5mm> 

\setplotarea x from -4 to 30, y from -3.4 to 2.0

\plot -4 -4       4 4         /         
\plot 4 4         12 -4          /      
\plot -2 -2       -2.95 -4 /            
\plot -.5 -2      -1.9 -4     /         
\plot -.5 -2      -.85 -4   /           
\plot 1 -2        .2 -4      /          
\plot 1 -2        1.25  -4     /        
\plot 2.5 -2      2.3  -4     /         
\plot 2.5 -2      3.35 -4      /        
\plot 5.5 -2      4.65  -4    /         
\plot 5.5  -2     5.7  -4      /        
\plot 7  -2       6.75   -4   /         
\plot 7 -2        7.8  -4      /        
\plot 8.5  -2     8.85  -4    /         
\plot 8.5  -2     9.9  -4      /        
\plot 10  -2      10.95  -4   /         
\plot 0  0        -.5  -2      /        
\plot 2 0         1 -2     /            
\plot 2 0         2.5   -2     /        
\plot 6 0         5.5    -2    /        
\plot 6 0         7 -2         /        
\plot 8 0         8.5 -2       /        
\plot 2 2         2 0         /         
\plot 6 2         6 0         /         

\multiput {\scriptsize $\bullet$} at   4 4   
       -4 -4    -2.95 -4    -1.9 -4    -.85 -4    .2 -4     1.25   -4   2.3  -4     3.35 -4    
        4.65  -4    5.7  -4     6.75  -4    7.8  -4      8.85  -4    9.9  -4      10.95  -4   12 -4     /

\put{$a_1$}[b] at 4 4.5
\put{$\bd^*S_1$}[t] at 4 -5


\plot 14  4       22 -4       /         
\plot 22 -4         30 4         /      
\plot 16 2       15.05 4 /              
\plot 17.5 2     16.1  4     /          
\plot 17.5 2     17.15  4   /           
\plot 19  2       18.2  4      /        
\plot 19 2         19.25   4     /      
\plot 20.5 2       20.3   4     /       
\plot 20.5  2      21.35 4      /       
\plot 23.5 2       22.65  4    /        
\plot 23.5  2      23.7   4      /      
\plot 25   2       24.75   4   /        
\plot 25  2        25.8  4      /       
\plot 26.5   2     26.85  4    /        
\plot 26.5   2     27.9   4      /      
\plot 28   2      28.95   4   /         
\plot 18  0        17.5  2      /       
\plot 20 0         19  2     /          
\plot 20 0         20.5   2     /       
\plot 24 0         23.5    2    /       
\plot 24 0         25 2         /       
\plot 26 0         26.5 2       /       
\plot 20 -2        20 0         /       
\plot 24 -2        24 0         /       

\multiput {\scriptsize $\bullet$} at   22 -4       
      14  4     15.05 4   16.1  4     17.15  4    18.2  4     19.25   4   20.3   4     21.35 4   
      22.65  4    23.7   4     24.75  4    25.8  4    26.85  4    27.9   4    28.95   4   30 4         /     

\put{$a_2$}[t] at 22  -4.5
\put{$\bd^*S_2$}[b] at 22  5

\put {$\circ$} at 8 0
\put {$\circ$} at 18 0
\plot 8.25 0  12.1 0 /
\plot 13.9 0 17.78 0 /
\plot    12.1   0    12.25 .4    12.25 -.4   12.55 .4   12.55 -.4
        12.85 .4   12.85 -.4  13.15 .4   13.15 -.4  13.45 .4
        13.45 -.4  13.75 .4   13.75 -.4  13.9 0     13.9  0 /

\endpicture
$$
\begin{center}\ \\
\emph{Figure 3}
\end{center}

We shall only be interested in computations on
tetrahedra with height $n \ge 2$.

If for $i=1,2$ we write $S^i = S^i_n(a_i) = \{ x_i \in \T^i : a_i \lle x_i\,,\;
d(x_i,a_i) \le n \}$, then 
$S = \{ x_1x_2 \in S_1 \times S_2 : \hor(x_1)+\hor(x_2)=0 \}$. 
The \emph{boundary} of $S^i$ in $\T^i$ is $\{a_i\} \cup \bd^* S^i$,
where $\bd^* S^1 = \{ b_1 \in S^1: b_1a_2 \in S \}$, respectively
$\bd^* S^2 = \{ b_2 \in S^2: a_1b_2 \in S \}$, 
and the boundary of $S$ in $\DL$ is
$$
\bd S= \bigl(\{ a_1 \} \times \bd^* S^2\bigr) \cup 
       \bigl(\bd^* S^1 \times \{ a_2 \}\bigr)\,.
$$
Imagining $S$ as a tetrahedron, two of its faces are copies of $S^1$ that meet
at the common bottom side $\bd^*S_1 \times \{a_2\}$, and the other two faces 
are  copies of $S^2$ that meet at the common top side $\{a_1\}\times \bd^*S_2$.
For $k=0,\dots, n$, the $k$-th \emph{level} of $S = S(a_1,a_2)$ is the
set
$$
L_k = L_k(a_1,a_2) = \{ x_1x_2 \in S : \hor(x_1) = \hor(a_1)+k \} 
= \bigl(S^1 \cap H^1_{\hor(a_1)+k}\bigr) \times 
  \bigl(S^2 \cap H^2_{\hor(a_2)+n-k}\bigr) \,.
$$
It has $q^k r^{n-k}$ elements. Furthermore, we write
$v_{1,s}$ for the successor of $a_1$ where the edge $[a_1,v_{1,s}]$ of
$\T^1$ has label $s \in \Z_q$, and analogously $v_{2,t}$ for the successor 
of $a_2$ where the edge $[a_2,v_{2,t}]$ of $\T^2$ has label $t \in \Z_r$. 

We identify functions on $S$ with their extensions to $\DL$, where the
latter have value $0$ on $\DL \setminus S$. In particular, every horizontal
function on $S$ must be $0$ on $\bd S$. We now choose non-zero functions
$\varphi^1$ on $\Z_q$ and $\varphi^2$ on $\Z_r$ such that
\begin{equation}\label{varphi}
\sum_{s=0}^{q-1} \varphi^1(s) = \sum_{t=0}^{r-1} \varphi^2(t) = 0
\AND
\sum_{s=0}^{q-1} \bigl(\varphi^1(s)\bigr)^2 
= \sum_{t=0}^{r-1} \bigl(\varphi^2(t)\bigr)^2 = 1\,.
\end{equation}
(Later on, we shall make specific choices for $\varphi^1$ and $\varphi^2$.)
Using these two functions, we define functions $f^1_k = f^1_k[\varphi^1]$
on $\T^1(a_1) = \{ x_1 \in \T^1 : a_1 \lle x_1 \}$ and 
$f^2_k = f^2_k[\varphi^2]$ on $\T^2(a_2) = \{ x_2 \in \T^1 : a_2 \lle x_2 \}$ by
$f^1_0 \equiv 0$, respectively $f^2_0 \equiv 0$, and for $k \ge 1$, 
\begin{equation}\label{f1f2}
\begin{aligned}
f^1_k(x_1) &= \begin{cases} \varphi^1(s)q^{(1-k)/2} 
                       &\text{if}\; v_{1,s} \lle x_1 \in H^1_{\hor(a_1)+k}\,,\\
                     0 & \text{otherwise\,,}
             \end{cases}
\AND\\[4pt]
f^2_k(x_2) &= \begin{cases} \varphi^2(t)r^{(1-k)/2} 
                       &\text{if}\; v_{2,t} \lle x_2 \in H^2_{\hor(a_2)+k}\,,\\
                     0 & \text{otherwise.}
             \end{cases}
\end{aligned}
\end{equation}
For $k \ge 1$, these functions have $\ell^2$-norm $1$. Now
we define for $0 \le k \le n$
\begin{equation}\label{fk}
f_{k,n}(x) = f_k[S,\varphi^1,\varphi^2](x) 
= f_k^1(x_1)f_{n-k}^2(x_2)\,,\quad x = x_1x_2 \in S=S(a_1,a_2)\,.
\end{equation}
Recall that $n = -\hor(a_2) - \hor(a_1) \ge 2$ is the height of $S$.
The following is a straightforward exercise.

\begin{lem}\label{action}
The functions $f_{k,n}$, $k=1,\dots, n-1$, are horizontal and
orthonormal in $\ell^2(S)$. The SRW-operator $P$ satisfies
$$
Pf_{k,n} = \frac{\sqrt{qr}}{q+r} \bigl( f_{k-1,n} + f_{k+1,n} \bigr).
$$
\end{lem}

Thus, since $f_{k,0}=f_{k,n} = 0$, the action of $P$ on the linear
space spanned by $f_{k,n}$, $k=1,\dots,n-1$, is described by the 
$(n-1)\times(n-1)$ tridiagonal matrix
\begin{equation}\label{matrix}
M_{n-1} = \frac{\sqrt{qr}}{q+r} 
\begin{pmatrix} 0 & 1      &        &  \\
                1 & \ddots & \ddots &  \\
                  & \ \ddots & \ \ddots & 1\\ 
                  &        & \ 1      & 0
\end{pmatrix}
\end{equation}
Its eigenvalues $\la_{m,n}$ and associated orthonormal eigenvectors
$\psi_{m,n}$ (the latter written as functions on $\{1,\dots, n-1\}$)
are
\begin{equation}\label{eigen}
\begin{aligned}
\la_{m,n} &= \frac{2\sqrt{qr}}{q+r} \cos \frac{m}{n}\pi \AND\\
\psi_{m,n}(k) &= \sqrt{\frac{2}{n}} \sin \frac{km}{n}\pi\,,\quad
m,k = 1, \dots, n-1
\end{aligned}
\end{equation}

\begin{cor}\label{ortho1}
The functions $g_{m,n}$ on $S$, $m=1,\dots, n-1$, defined by
$$
g_{m,n} = \sum_{k=1}^{n-1} \psi_{m,n}(k)f_{k,n}
$$
are horizontal and orthonormal in $\ell^2(S)$ as well as in $\ell^2(DL)$.
They satisfy
$$
\begin{aligned}
Pg_{m,n} &= \la_{m,n} \cdot g_{m,n}\qquad\mbox{and}\\
\spn\{ g_{m,n} : m=1,\dots,n-1 \} &= \spn\{ f_{k,n} : k=1,\dots,n-1 \}\,.
\end{aligned}
$$
\end{cor}

Once more, recall that besides depending on the height $n$ of $S$, 
each $g_{m,n}$ depends on  $a_1$, $a_2$, $\varphi^1$ and $\varphi^2$,
$$
g_{m,n} = g_m[S,\varphi^1,\varphi^2]\,.
$$

\begin{lem}\label{ortho2}
Let $S(a_1,a_2)$ and $S(\tilde a_1,\tilde a_2)$ be two tetrahedra
of heights $n$ and $\tilde n \ge 2$, respectively. 
Let $\varphi^1$ and $\varphi^2$, respectively $\tilde\varphi^1$ and 
$\tilde\varphi^2$ be as in \eqref{varphi}. Write
$f_{k,n} = f_k[a_1,a_2,\varphi^1,\varphi^2]$ and $\tilde f_{l,\tilde n} 
= f_l[\tilde a_1,\tilde a_2,\tilde \varphi^1,\tilde \varphi^2]$.
\quad If one of 

\begin{itemize}

\item[(i)] $(a_1,a_2) \ne (\tilde a_1,\tilde a_2)\,,\quad$ or

\item[(ii)] $(a_1,a_2) = (\tilde a_1,\tilde a_2)\;\ $ and 
$\ \;\varphi^1 \perp \tilde\varphi^1$ or $\varphi^2 \perp \tilde\varphi^2$
 
\end{itemize}
holds, then
$$
\spn\{ f_{k,n} : k=1,\dots,n-1 \} \perp 
\spn\{ \tilde f_{l,\tilde n} : l=1,\dots,\tilde n-1 \}\,.
$$
\end{lem}


\begin{proof} (i) If $S(a_1,a_2) \cap S(\tilde a_1,\tilde a_2) =
\emptyset$ then perpendicularity is obvious.

If $S(a_1,a_2) \cap S(\tilde a_1,\tilde a_2) \ne \emptyset$ then 
both $a_1,\tilde a_1$ and $a_2,\tilde a_2$
must be comparable with respect to the partial order $\lle$.  
Assume that $a_1 \lle \tilde a_1$ and $a_1 \ne \tilde a_1$. (The other
three cases are treated analogously). Let $\ka = \hor(\tilde a_1) - \hor(a_1)$.

If $k \ne \ka + l$ then we certainly have $f_{k,n} \perp \tilde f_{l,\tilde n}$
(since the two functions have disjoint support).

If $k = \ka+ l$ then by construction, the function $f^1_k$ on $\T^1(a_1)$,
as defined in \eqref{f1f2}, is constant on $\supp \tilde f^1_l$,
where the latter is given as in \eqref{f1f2}, but on $\T^1(\tilde a_1)$.
Since $\sum \tilde f^1_l = 0$, we have $f^1_k \perp \tilde f^1{l}$,
Therefore, with $f^2_k$ and $\tilde f^2_l$ defined in the same way on 
$\T^2(a_2)$, respectively $\T^2(\tilde a_2)$,
$$
\sum_x f_{k,n}(x)\tilde f_{l,\tilde n}(x) 
= \sum_{x_2} f^2_{n-k}(x_2) \tilde f^2_{\tilde n - l}(x_2) 
\underbrace{\sum_{x_1} f^1_k(x_1) \tilde f^2_l(x_1)}_{\dps =0} = 0\,.
$$
(ii) If $(a_1,a_2) = (\tilde a_1,\tilde a_2)$ then $\tilde n = n$, 
and $f_{k,n} \perp \tilde f_{l,n}$ if $l \ne k$, since the
functions have disjoint support. Finally,
$$
\hspace*{4.6cm} \langle f_{k,n},  \tilde f_{k,n} \rangle 
= \langle \varphi^1 , \tilde \varphi^1 \rangle
\langle \varphi^2 , \tilde \varphi^2 \rangle = 0 \hspace{4.6cm}\qedhere
$$
\end{proof}

We now specify our choices for the functions $\varphi^1$ and $\varphi^2$.
For $i = 1, \dots, q-1$ and $s \in \Z_q$, let
\begin{equation}\label{phii}
\varphi_i^1(s) = 
\begin{cases} \dfrac{q-i}{\sqrt{(q-i)(q+1-i)}}\,,& s=i-1\,,\\
              -\dfrac{1}{\sqrt{(q-i)(q+1-i)}}\,,& s=i, \dots, q-1\,,\\
              0\,,& \mbox{otherwise}.
\end{cases}
\end{equation}
These functions are orthogonal and satisfy \eqref{varphi}.
Analogously, replacing $i$ with $j$ and $q$ with $r$, we define the
orthogonal functions $\varphi_j^2(t)$, $j=1, \dots, r-1$ ($t \in \Z_r$).
We shall write
\begin{equation}\label{fkgm}
f_k^{[S,i,j]} = f_k[S,\varphi_i^1,\varphi_j^2] \AND
g_m^{[S,i,j]} = g_m[S,\varphi_i^1,\varphi_j^2]
\end{equation}

\begin{pro}\label{basis}
The set 
$$
{\mathfrak B}_S =\Bigl\{ g_m^{[\tilde S,i,j]} \ :\ \tilde S \subseteq S\,;\; 
m\in\{1,\dots,n(\tilde S) - 1\}\,;\; 
i\in\{1,\dots, q-1\}\,;\;j\in\{1,\dots,r-1\} \Bigr\}
$$
constitutes an orthonormal basis of the linear space of all horizontal
functions on the tetrahedron $S$. Here, $\tilde S$ runs through all
tetrahedra in $\DL(q,r)$ that are contained in $S$ and have height
$n(\tilde S) \ge 2$.
\end{pro}

\begin{proof} Instead of the functions $g_m^{[\tilde S, i,j]}$, we may equivalently
work with the functions $f_m^{[\tilde S,i,j]}$, since they are also linearly
independent and span the same space as ${\mathfrak B}_S$.

If $L_k$ is the $k$-th level of $S$ (where $S$ is assumed to have height $n$) and
$f$ is any horizontal function with support in $L_k$, 
then $f$ must satisfy each of the following $q^k + r^{n-k}$ equations
$$
\sum_{x_2:x_1x_2 \in L_k} f(x_1x_2) = 0  \AND
\sum_{x_1:x_1x_2 \in L_k} f(x_1x_2) = 0\,.
$$
Thus, the dimension of the linear space of all horizontal functions with 
support in $L_k$ is $(q^k-1)(r^{n-k}-1)$. 

On the other hand, we can count all $f_m^{[\tilde S, i,j]}$
where $\tilde S \subseteq S$ and $\supp f_m^{[\tilde S, i,j]} \subseteq L_k$.
We find $(q-1)(r-1)$ functions of this type (one for each pair $(i,j)$)
associated with every tetrahedron $\tilde S = S(\tilde a_1,\tilde a_2)$,
where $a_1 \lle \tilde a_1$ and $d(\tilde a_1, a_1) \le k-1$, and
at the same time $a_2 \lle \tilde a_2$ and $d(\tilde a_2, a_2) \le n-k-1$.
There are precisely
$$
\sum_{\ka=0}^{k-1} q^{\ka} \sum_{\nu=0}^{n-k-1} r^{\nu} 
= \frac{(q^k-1)(r^{n-k}-1)}{(q-1)(r-1)} 
$$
choices for $(\tilde a_1,\tilde a_2)$. Thus, the number of all functions
$f_m^{[\tilde S,i, j]}$ with support in $L_k$ (which are linearly
independent) coincides with the dimension of the  space of all horizontal 
functions on $L_k$, and we have a basis of that space.
Putting together the different levels of $S$, we obtain the proposed 
result.
\end{proof}

Thus, we obtain the following spectral decomposition of $\ell^2(\DL)$.

\begin{thm}\label{spec} The spectrum of of the SRW-operator $P$ on
$\DL(q,r)$ is the interval $[-\rho(P)\,,\,\rho(P)]$, where
$\rho(P) = 2\sqrt{qr}/(q+r)$.

It is a pure point spectrum, being the closure of the set 
$$
\left\{\la_{m,n} = \rho(P) \cos \frac{m}{n}\pi : n \ge 2\,,\;1 \le m \le n-1
\right\}.
$$

Furthermore, the set of all functions
$$
{\mathfrak B} = \left\{ g_m^{[S,i,j]} \ :\ 
{  S\;\mbox{tetrahedron in}\; \DL\,;\;
m\in\{1,\dots,n(\tilde S) - 1\}\,;\; \atop
i\in\{1,\dots, q-1\}\,;\;j\in\{1,\dots,r-1\}  } \right\}\,,
$$
constructed in Corollary \ref{ortho1}, respectively \eqref{fkgm},
is a complete orthonormal system in $\ell^2(\DL)$ consisting of finitely
supported functions; we have
$$
Pg_m^{[S,i,j]} = \la_{m,n} \cdot g_m^{[S,i,j]}\,,
$$
where $n = n(S)$ is the height of $S$.
\end{thm}

\begin{proof}
If $f$ is any horizontal function on $\DL$, then there is some
tetrahedron $S$ containing its support. By Proposition \ref{basis},
$f$ is a linear combination of elements of ${\mathfrak B}_S$. 
Thus, ${\mathfrak B} = \bigcup_S {\mathfrak B}_S$ is an orthonormal
system that spans the space of all horizontal functions.
Now Lemma \ref{dense} completes the proof.
\end{proof}

\vspace{.4cm}
\medskip

\section{The spectral measures}\label{meas}

Using Theorem \ref{spec}, we can compute the spectral measures
\eqref{measures}. Indeed, if $x \in \DL(q,r)$, then the Fourier 
expansion of $\de_x$ with respect to the orthonormal system
of \eqref{spec} is
$$
\de_x = \sum_{S} \sum_{m=1}^{n(S)-1} \sum_{i=1}^{q-1} \sum_{j=1}^{r-1} 
        g_m^{[S,i,j]}(x)\, g_m^{[S,i,j]}(\cdot)\,.
$$
Therefore, for $x,y \in \DL$,
$$
p^{(N)}(x,y) = \langle \de_x,P^N\de_y \rangle 
= \sum_{S} \sum_{m=1}^{n(S)-1} \sum_{i=1}^{q-1} \sum_{j=1}^{r-1} 
\la_{m,n(S)}^N\, g_m^{[S,i,j]}(x)\, g_m^{[S,i,j]}(y)\,,
$$
and comparing this with \eqref{moments}, we find that for any continuous
function $\ff$ on $\spec(P)$, its integral with respect to $\mu_{x,y}$
is
\begin{equation}\label{muxy}
\int_{\spec(P)} \ff(\la)\,\mu_{x,y}(d\la) = 
\sum_{S} \sum_{m=1}^{n(S)-1} \sum_{i=1}^{q-1} \sum_{j=1}^{r-1} 
\ff(\la_{m,n(S)})\, g_m^{[S,i,j]}(x)\, g_m^{[S,i,j]}(y)\,,
\end{equation}
a countable sum of point masses. Since $\DL$ is transitive,
we only need the measures $\mu_{o,x}$, where $x \in \DL$.
Let $x=x_1x_2$ with $\up_i = \up(x_i)$ and $\dn_i=\dn(x_i)$.
Furthermore, for $i=1,2$ let $c_i = x_i \cf o_i$, so that
$S(c_1,c_2)$ has height $\nn = \up_1+\up_2 = \dn_1 + \dn_2 = 
\bigl(d(x_1,o_1)+d(x_2,o_2)\bigr)/2$, see Figure 4.

$$
\beginpicture
\setcoordinatesystem units <3.2mm,3.7mm> 

\setplotarea x from -4 to 24, y from -3.8 to 6.4

\put{$o_1$} [lb] at  8.15 0.2
\put{$x_1$} [rt] at 1.85 1.85
\put{$c_1$} [lb] at 4.15 4.2 
\put{$a_1$} [lb] at 2.15 6.2 


\plot 8 0  4 4  2 2 / 
\plot 4 4  2 6 /

\put{$o_2$} [rt] at  11.95 -.2
\put{$x_2$} [lb] at  22.15  2.2
\put{$c_2$} [rt] at  15.95 -4.2
\put{$a_2$} [rt] at  17.95 -6.2


\plot 12 0   16 -4   22 2 /         
\plot 16 -4  18 -6 /

\multiput {$\scs\bullet$} at 8 0  12 0  2 2  22 2  4 4  16 -4  
                             2 6  18 -6  /

\arrow <5pt> [.2,.67] from 0.4 7.6 to 0 8
\arrow <5pt> [.2,.67] from 19.6 -7.6 to 20 -8
\setdashes <1.5pt>
\plot  0.4 7.6    2 6 /
\plot 19.6 -7.6  18 -6 /
\put{$\omega_1$} [rb] at 0.2 8.2
\put{$\omega_2$} [lt] at 20.2 -8.2


\setdots<4pt>
\plot -4.5 0  12 0 /
\plot -1.5 2  22 2 /
\plot -4.5 4  4 4 /
\plot -4.5 -4  16 -4 /

\put{$\up_1\!\!\left\{\rule[-6mm]{0mm}{0mm}\right.$}[rb] at -4.5 -.1 
\put{$\up_2\!\!\left\{\rule[9mm]{0mm}{0mm}\right.$}[rt] at -4.5 0.1
\put{$\dn_1\!\!\left\{\rule[5mm]{0mm}{0mm}\right.$}[rb] at -1.5 2
\put{$\dn_2\!\!\left\{\rule[-9.95mm]{0mm}{0mm}\right.$}[rt] at -1.5 2.1

\endpicture
$$
\vspace{.1cm}

\begin{center}
\emph{Figure 4}
\end{center}

\vspace{.4cm}

Since $\la_{m,n}$ depends only on $m/n$, we choose
$m$ and $n$ relatively prime ($1 \le m < n$). 
In order to compute the mass $\mu_{o,x}(\la_{m,n})$, we have to consider
all tetrahedra $S=S(a_1,a_2)$ with height $n(S) = \ell n$ that contain 
$S(c_1,c_2)$, that is, $a_i \in \geo{c_i,\om_i}$ for $i=1,2$. For those $S$, 
note that $g_m^{[S,i,j]}(o) = 0$ when $(i,j) \ne (1,1)$. 
Therefore, denoting by $\lceil\cdot\rceil$ the next larger integer, 
\eqref{muxy} now yields
\begin{equation}\label{muox1}
\mu_{o,x}(\la_{m,n}) = \sum_{\ell=\lceil\nn/n\rceil}^{\infty} 
\sum_{\scs S \supseteq S(c_1,c_2) \atop \scs n(S) = \ell n} 
g_{\ell m}^{[S,1,1]}(o)\, g_{\ell m}^{[S,1,1]}(x)\,,
\end{equation}
and if we set 
$$
k_i = k_i(S) = d(a_i,c_i) = \hor(c_i)-\hor(a_i)\,,\quad  i=1,2\,,
$$
then by Corollary \ref{ortho1}
\begin{equation}\label{muox2}
\begin{aligned}
g_{\ell m}^{[S,1,1]}(x) \,g_{\ell m}^{[S,1,1]}(o) =&\  
C_{k_1,k_2} \,\, q^{-k_1-(\up_1+\dn_1)/2} \, r^{-k_2-(\up_2+\dn_2)/2} \;\times\\
&\times \,\frac{2}{\ell n}\, \sin \Bigl((k_1+\up_1)\tfrac{m}{n}\pi\Bigr) \,\,
\sin \Bigl((k_1+\dn_1)\tfrac{m}{n}\pi\Bigr)\,.
\end{aligned}
\end{equation}
where $n(S) = \ell n$ and
\begin{equation}\label{muox3}
C_{k_1,k_2} = 
\begin{cases}
(q-1)(r-1)
\,, &\mbox{if $\;k_1 > 0\;$ and $\;k_2 > 0\,$,}\\
-(r-1)
\,, &\mbox{if $\;k_1 = 0\;$ and $\;k_2 > 0\,$,}\\
-(q-1)
\,, &\mbox{if $\;k_1 > 0\;$ and $\;k_2 = 0\,$,}\\
1
\,, &\mbox{if $\;k_1 = k_2 = 0\,$.}
\end{cases}
\end{equation}
The last case occurs only when $a_1=c_1$ and $a_2=c_2$.
Also, note that $k_1 + k_2 = n(S) - \nn$ and that $\up_i+\dn_i = d(x_i,o_i)$,
$i = 1,2$. We obtain

\begin{pro}\label{muox} If $1 \le m < n$ and $m$ and $n$ are 
relatively prime, then with constants as in \eqref{muox3},
$$
\begin{aligned}
&\mu_{o,x}(\la_{m,n})\\ &=\!\! \sum_{\ell=\lceil\nn/n\rceil}^{\infty}
\frac{2r^{-\ell n}}{\ell n}\, \sum_{k=0}^{\ell n -\nn} C_{k, \ell n- \nn- k}\,
(r/q)^{k+d(x_1,o_1)/2}\,\sin \Bigl((k+\up_1)\tfrac{m}{n}\pi\Bigr) \,
\sin \Bigl((k+\dn_1)\tfrac{m}{n}\pi\Bigr)\,.
\end{aligned}
$$
\end{pro}

Elementary computations yield the following.

\begin{cor}\label{muoo}
The Plancherel measure $\mu = \mu_{o,o}$ is given by
$$
\mu(\la_{m,n}) = \begin{cases}
\dfrac
{\log(1-r^{-n}) - \log(1-q^{-n})}{n(r-q)}\,\,
\dfrac{2qr(q+r)(q-1)(r-1) \sin^2 \frac{m}{n}\pi}
{(r-q)^2 + 4qr \sin^2 \frac{m}{n}\pi}
                 \,,&\mbox{if $\;r\ne q\,$,}\\
\dfrac{(q-1)^2}{q^n-1}\,,&\mbox{if $\;r = q\,$,}                 
\end{cases}
$$
where $\la_{m,n} = \frac{2\sqrt{qr}}{q+r} \cos \frac{m}{n}\pi$ and
$m$ and $n$ are relatively prime ($1 \le m < n$). 
\end{cor}

Note that $\la_{n-m,n} = -\la_{m,n}$ and $\mu(\la_{n-m,n}) = \mu(-\la_{m,n})$, 
that is, the Plancherel measure is symmetric, as it has to be, since 
$\DL(q,r)$ is a bipartite graph (it has no odd cycles).
The formula for $\mu$ in the case $r=q$ (lamplighter group) was obtained 
previously in \cite{GrZu2} and \cite{DiSc}. 
For $x \ne o$,
the inner sum in Proposition \ref{muox} can be computed in (lengthy)
closed form, but in general not the outer one.
 
\section{Asymptotic behaviour of the return probabilities}\label{asymptotics}

Combining \eqref{moments} with Corollary \ref{muoo}, we can determine the
exact asymptotic behaviour of the return probabilities $p^{(N)}(o,o)$
as $N \to \infty$. For odd $N$, these probabilities are $=0$.

For dealing with $p^{(2N)}(o,o)$, the following standard technical lemma will
be useful.

\begin{lem}\label{technical} For $k \in \N$ and $\ga \in \R$ and any
sequence $\ep_n$ tending to\/ $0$, let
$$
\Sig(N) = \Sig(N;k,\ga) = 
\sum_{n=2k+1}^{\infty} (1+\ep_n)\, n^{\ga} \,q^{-n}\, \cos^{2N}
\tfrac{k}{n}\pi\,.
$$
Then
$$
\Sig(N) \sim 
\xi_k^{\ga}\,(2\pi/C_k)^{1/2} \, \exp\bigl(-B_k\,N^{1/3}\bigr) \, N^{(1+2\ga)/6} \quad
\mbox{as}\;N \to \infty\,,
$$
where
$$
\begin{gathered}
\Phi_k(\xi) = \xi \log q + \frac{(k\pi)^2}{\xi^2},\quad
\xi_k = \left( \frac{2(k\pi)^2}{\log q} \right)^{1/3};\\
B_k = \Phi_k(\xi_k) = 3 \left( \frac{k\pi\log q}{2} \right)^{2/3},\quad
\Phi'_k(\xi_k)=0 \AND
C_k = \Phi_k''(\xi_k) = 6 \left( \frac{(\log q)^2}{4k\pi} \right)^{2/3}.
\end{gathered}
$$
\end{lem}

(Here, as usual,  $\sim$ denotes asymptotic equivalence,  
i.e., quotients tending to $1$. It will always be clear from
the context whether $\sim$ means neighbourhood in a graph or asymptotic
equivalence.)

\begin{proof}
We decompose $\Sig(N) = \Sig_1(N) + \Sig_2(N)$, where the sum
$\Sig_1(N)$ ranges over all $n$ with $2k < n < \ka N^{1/3}$, and
$\Sig_2(N)$ ranges over all $n \ge \ka N^{1/3}$, and where
$\ka = \xi_k/2$. 

We start with $\Sig_2(N)$. For $n \to \infty$,
$$
\cos^{2N} \tfrac{k}{n}\pi  = 
\exp\Bigl(-(k\pi)^2 N\bigl(n^{-2} + O(n^{-4})\bigr)\Bigr)\,.
$$
Therefore, we have as $N \to \infty$
\begin{equation}\label{Sig2first}
\begin{aligned}
\Sig_2(N) &\sim N^{\ga/3} \sum_{n \ge \ka N^{1/3}} 
\xi_{N,n}^{\ga} 
\exp \left( -N^{1/3} \Phi_k (\xi_{N,n})\right)\,, 
\quad\mbox{where}\quad\xi_{N,n} = \frac{n}{N^{1/3}}.
\end{aligned}
\end{equation}
The point where $\Phi_k$ attains its minimum is $\xi_k$, and we
compute the values $\Phi_k(\xi_k) = B_k$ and $\Phi_k''(\xi_k) = C_k$,
as given above.  Therefore
$$
\Phi_k(\xi) = B_k + \frac{C_k}{2}(\xi - \xi_k)^2 + R(\xi)\,(\xi - \xi_k)^3
$$
with $R(\xi)$ continuous for $\xi > 0$. We now ``substitute''
$$
\tau_{N,n} = N^{1/6}(\xi_{N,n}-\xi_k)\,, \quad\mbox{with}\quad 
\Delta \tau_{N,n} = \tau_N(n+1) - \tau_{N,n} = N^{-1/6} \to 0\,.
$$
Then we can rewrite \eqref{Sig2first} as
\begin{equation}\label{Sig2second}
\begin{aligned}
\Sig_2(N) \sim&\ \exp\bigl(-B_k\,N^{1/3}\bigr)\, N^{(1+2\ga)/6}\,\times\\
&\times\!\! \sum_{\scs n :\atop \scs \tau_{N,n} \ge -\ka N^{1/6}}
(\xi_k+ N^{-1/6}\tau_{N,n})^{\ga} \,
\exp \left\{ -\tfrac{C}{2}\tau_{N,n}^2 - R^*_N(\tau_{N,n}) \right\}\,
\Delta \tau_{N,n}\end{aligned}
\end{equation}\\[1pt]
where $R^*_N(\tau) = N^{-1/6} R(\xi_0+ N^{-1/6}\tau)\, \tau^3$.

It is now standard that the sum in \eqref{Sig2second}
converges to 
$\xi_k^{\ga} \int_{-\infty}^{\infty} e^{-C_k\tau^2/2}\,d\tau 
= \xi_k^{\ga}\sqrt{2\pi/C}$.
(One has to use dominated convergence in a suitable central piece of the
sum and control the two tails.)
Thus, $\Sig_2(N)$ has the asymptotic behaviour that we have proposed for
$\Sig(N)$.

\vspace{.3cm}

Let us now look at $\Sig_1(N)$. Set 
$M = \sup_n |1+\ep_n|\, n^{\ga}\, q^{-n}$. Then
$$
\Sig_1(N) \le M\, \ka \, N^{1/3} \, \cos^{2N} 
\left(\frac{k\pi}{\ka N^{1/3}}\right) 
\sim M\, \ka \, N^{1/3} \, \exp \left(-\frac{(k\pi)^2}{\ka^2} N^{1/3} \right)\,.
$$
With our choice $\ka = \xi_k/2$, one checks that $(k\pi)^2/\ka^2 > B_k$.
Therefore $\Sig_1(N)/\Sig_2(N) \to 0$ as $N \to \infty$.
\end{proof}

For the following, recall from Theorem \ref{spec} that $\rho(P) = 1$
when $r=q$. When $r \ne q$, it is enough to consider only $r > q$.

\begin{thm}\label{return} Let $\xi_1, B_1$ and $C_1$ be as defined in Lemma
\ref{technical}.\\[3pt]
{\rm (i)} If $r > q$ then
$$
p^{(2N)}(o,o) \sim 
A_1\,\rho(P)^{2N} \,
\exp\bigl(-B_1\,N^{1/3}\bigr)\, N^{-5/6}
\qquad{as}\;\; N \to \infty\,, 
$$
where
$$
A_1 = 4\,\pi^2\,\xi_1^{-3}\, (2\pi/C_1)^{1/2}\,qr(q+r)(q-1)(r-1)\big/(r-q)^3\,.
$$\\[3pt]
{\rm (ii)} If $r=q$ then 
$$
p^{(2N)}(o,o) \sim 
\bar A_1 \, \exp\bigl(-B_1\,N^{1/3}\bigr)\, N^{1/6}
\qquad{as}\;\; N \to \infty\,, 
$$
where
$$
\bar A_1 = 2(q-1)^2\,(2\pi/C_1)^{1/2}\,.
$$
\end{thm}

\begin{proof}
We decompose (using $\la_{n-m,n} = -\la_{m,n}$ and $\la_{1,2}=0$)
$$
p^{(2N)}(o,o) = 
\sum_{n=3}^{\infty} \sum_{m=1 \atop\scs  \gcd(m,n)=1}^{n-1}
\mu(\la_{m,n})\, \la_{m,n}^{2N} = S_1(N) + S_2(N)\,,
$$
where
$$
S_1(N) = 2\sum_{n=3}^{\infty} \mu(\la_{1,n})\, \la_{1,n}^{2N}
\AND
S_2(N) = \sum_{n=4}^{\infty} \sum_{\scs m=2 \atop\scs  \gcd(m,n)=1}^{n-2}
\mu(\la_{m,n})\, \la_{m,n}^{2N}\,.
$$\\[3pt]
\underline{Case $r > q$.} 
Then, for $n \to \infty$,
$$
2 \,\mu(\la_{1,n}) = 2\,\pi^2\, A_0 \,n^{-3}\,q^{-n} \bigl(1+O(n^{-2})\bigr)
\quad\mbox{where}\quad
A_0 = 2 qr(q+r)(q-1)(r-1)\big/(r-q)^3\,.
$$
Therefore, using Lemma \ref{technical}, 
\begin{equation}\label{S1}
S_1(N) = A_0\,\rho(P)^{2N} \,\Sig(N;1,-3)
\end{equation}
has the asymptotic behaviour that we have proposed for $p^{(2N)}(o,o)$.
Thus, it remains to show that $S_2(N)$ is dominated by $S_1(N)$.
Note that for $2 \le m \le n-2$,  we have 
$$
\la_{m,n}^2 \le \rho(P)^2 \,\cos^2 \tfrac{2}{n}\pi\,, 
$$
and 
$$
\sum_{\scs m=2 \atop\scs  \gcd(m,n)=1}^{n-2} \mu(\la_{m,n}) \le 
A_0 \bigl(\log(1-r^{-n}) - \log(1-q^{-n})\bigr)
= A_0 \,(1+ \ep_n)\, q^{-n}\,,
$$
where $\ep_n \to 0$ as $n \to \infty\,$. Therefore,
\begin{equation}\label{S2}
S_2(N) \le A_0 \,\rho(P)^{2N} \,\Sig(N;2,0)\,.
\end{equation}
Since $B_2 > B_1$, Lemma \ref{technical} shows that 
$\Sig(N;2,0)/\Sig(N;1,-3) \to 0$, and comparing \eqref{S1} with \eqref{S2},
we see that $S_2(N)/S_1(N) \to 0$ as $N \to \infty$.\\[3pt]
\underline{Case $r = q$.} The proof is basically the same. The only
difference is that in this case, Theorem \ref{spec}  yields
$$
\mu(\la_{1,n}) \sim (q-1)^2 \,q^{-n}\,,
$$
while in case $r > q$, we had the additional factor $n^{-3}$. Therefore
$$
S_1(N) = 2(q-1)^2 \,\Sig(N;1,0) \AND S_2(N) \le (q-1)^2 \, \Sig(N;2,1)\,,
$$
whence Lemma \ref{technical} implies the result.
\end{proof}

It is quite surprising that the polynomial terms $N^{-5/6}$ versus
$N^{1/6}$ are different when $r \ne q$, respectively $r=q$.
The asymptotics in case $r=q$ were computed previously by {\sc Revelle} 
\cite{Re} for a very similar random walk on the lamplighter group: in terms of
$\DL(q,r)$ this is SRW on the graph obtained by adding edges in the first tree
($\T_q$), so that each vertex is connected to each of the siblings of its
predecessor and then taking the horocyclic product as before. It turns out
that for $q=r$, the $N$-step return probabilities of that random walk
are just $2 \,p^{(N)}(o,o)$, see \cite{Wo}, \S 5 for details. 
The specific computations of \cite{Re} are very similar to ours, although 
\cite{Re} does not use the Plancherel measure.

\section{Plancherel measure and cumulative spectral measure}\label{cumulative}

Let $(X,o)$ and $(X',o')$ be locally finite, infinite graphs with respective
roots $o$ and $o'$. Let $R(X,X')$ be the largest radius $R$ for which there is
a root-preserving isomorphism between the balls $B_X(o,R)$ and 
$B_{X'}(o',R)$ in the respective graph metrics.
A sequence of rooted graphs $(X_n,o_n)$ is said to \emph{converge} to
the graph $(X,o)$, if $R(X,X_n) \to \infty$. 
In this situation, given vertices $x, y \in X$, we can consider them 
via the respective isomorphisms as elements of $X_n$ for all $n \ge n(x,y)$,
and we have the corresponding spectral measures $\mu_{x,y}^{(n)}$ associated
with the SRW operator $P_n$ on $X_n$, as well as the measure $\mu_{x,y}$
associated with $P$ on $X$. (In particular, we consider $o$ as the common 
root of all graphs in the sequence.)

In this setting, it is a well known fact regarding operator convergence that
\begin{equation}\label{weak}
\mu^{(n)}_{x,y} \to \mu_{x,y} \quad\mbox{weakly, as}\; n \to \infty\,;
\end{equation}
see e.g. {\sc Grigorchuk and \.Zuk} \cite{GrZu1}. The main interest
here is in the diagonal elements $\mu^{(n)}_{x,x} \to \mu_{x,x}$,
in particular in the case when $X$ is vertex-transitive, and
$\mu_{x,x} = \mu_{o,o} \;\forall\ x \in X$.

When the $X_n$ are \emph{finite} graphs, another type of spectral measure 
is of interest in the place  of $\mu^{(n)}_{o,o}$, namely the
\emph{cumulative spectral measure}
$$
\tilde \mu_{X_n} = \tilde \mu_n =\frac{1}{|X_n|}\sum_{x \in X_n} \mu^{(n)}_{x,x} 
= \frac{1}{|X_n|}\sum_{\la \in \spec(P_n)} \mlt(\lambda)\, \de_{\lambda}\,,
$$
where $\mlt(\lambda)$ is the multiplicity of $\lambda$ as
an eigenvalue of $P_n$. If $X_n$ is vertex-transitive then $\tilde \mu_n = 
\mu^{(n)}_{o,o}$. In general, the following two questions are of interest.
\begin{equation}\label{question}
\begin{aligned}
&\mbox{(a) Does the sequence $(\tilde \mu_n)$ converge weakly to some 
probability measure $\tilde \mu$ ?}\\
&\mbox{(b) Is $\tilde \mu = \mu_{o,o}$~?}
\end{aligned}
\end{equation}
In case of a positive answer to question (\ref{question}.a),
we call $\tilde \mu$ the \emph{cumulative spectral measure}
of $(X,o)$ with respect to the sequence $(X_n,o)$. Recently, the names
Von Neumann and Serre have been associated with that measure.

All this applies, in particular, when each $X_n$ is an induced subgraph of $X$.
(``Induced'' means that when $x,y \in X_n$ are neighbours in $X$, then 
also $x \sim y$ in $X_n$.) In this situation, it may have advantages
to use instead of $P_n$ the restriction (truncation) $P_{|n} = P|_{X_n}$ 
of $P$ on $X_n$. Note that $P_n$ and $P_{|n}$ coincide in the interior
of $X_n$, while they differ in the points of the boundary 
$\partial X_n$ of $X_n$ (i.e., the points of $X_n$
having a neighbour in $X \setminus X_n$), where $P_{|n}$ is strictly
substochastic. The operator $P_{|n}$ acts on the same $\ell^2$-space as
$P$, while $P_n$ uses different weights (vertex degrees) at the boundary
points. The  spectral measures of $P_{|n}$ also converge weakly to
the respective spectral measures of $X$, and one can as well study
the cumulative spectral measures $\tilde \mu_{|n}$ associated with 
$P_{|n}$ and their possible limit, again denoted $\tilde \mu$.

We remark that in the literature, the distinction between
Plancherel (Kesten) and cumulative spectral measures has not always been
very clear. 

Question (\ref{question}) has first been dealt with explicitly by
{\sc McKay} \cite{McK}, who showed that when the $X_n$ are $(q+1)$-regular
graphs with (asymptotically) few cycles, then the sequence
$\tilde\mu_n$ converges weakly to the Plancherel measure of the 
tree $\T_q$. A systematic answer to question (\ref{question}.a) is given
by {\sc Serre} \cite{Ser}. In \cite{GrZu2}, the spectrum
of SRW on the lamplighter group $\Z_2 \wr \Z$ (i.e., $DL(2,2)$)
is computed via an approximating sequence of Schreier graphs, and it is
shown that the corresponding cumulative spectral measure coincides
with $\mu_{o,o}$.

The sequence $(X_n)$ of subsets of $X$ is called a \emph{F\o lner
sequence,} if $\inf_n |\partial X_n|/|X_n| = 0$. It is called 
\emph{expanding,} if that infimum is positive. Recall that a graph is called
\emph{amenable} if it has a F\o lner sequence. By {\sc Dodziuk} \cite{Dod},
this holds if and only if the spectral radius of SRW satisfies $\rho(P)=1$.
This notion comes from group theory; a group is called 
amenable if it carries a finitely additive, left-invariant probability
measure, and a finitely generated group is amenable if and only if one
(equivalently, each) of its Cayley graphs with respect to a finite,
symmetric generating set is amenable. For more details in the context of
random walks, see e.g. \cite{Wbook}, \S 10 and \S 12.

\begin{thm}\label{expanding} Let $X$ be an infinite, connected 
locally finite vertex-transitive graph, and $(X_n)$ an increasing 
subsequence of finite subgraphs whose union is $X$.
\\[3pt]
{\rm (a)} If $(X_n)$ is a F\o lner sequence then both $\tilde \mu_n$
and $\tilde \mu_{|n}$ converge weakly to the Plancherel measure
of $X$.
\\[3pt]
{\rm (b)} If $(X_n)$ is an expanding sequence then $\tilde \mu_{|n}$ 
does \emph{not} converge to the Plancherel measure of $X$.
\end{thm}

\begin{proof} All involved measures are probability measures 
on $\R$ with compact support. Therefore, weak convergence holds
if and only if for each $N \in \N$, the $N$-th moment of $\tilde \mu_n$
(resp. $\tilde \mu_{|n}$) converges to the $N$-th moment of
the Plancherel measure. The latter is 
$M_N(\mu_{o,o})=p^{(N)}(o,o) =  p^{(N)}(x,x)$ for all $x \in X$
(by transitivity), while 
$$
M_N(\tilde\mu_n) = \frac{1}{|X_n|} \sum_{x \in X_n} p_{X_n}^{(N)}(x,x)
$$
(analogously for $\tilde\mu_{|n}$). Now consider
$$
\partial_N\! X_n = \{ x \in X_n : d(x,\partial X_n) \le N/2 \} 
$$
If $x \in X_n \setminus \partial_N\! X_n$ then
$p_{X_n}^{(N)}(x,x) = p_{|X_n}^{(N)}(x,x) = p^{(N)}(x,x)$.
Therefore
$$
\begin{aligned}
M_N(\mu_{o,o})- M_N(\tilde\mu_{|n}) &= \frac{1}{|X_n|} \sum_{x \in X_n} 
\Bigl( p^{(N)}(x,x) - p_{|X_n}^{(N)}(x,x) \Bigr) \\
&= \frac{1}{|X_n|} \sum_{x \in \partial_N\! X_n} 
\Bigl( p^{(N)}(x,x) - p_{|X_n}^{(N)}(x,x) \Bigr)
\end{aligned}
$$
Now, if $(X_n)$ is a F\o lner sequence then 
$$
|M_N(\mu_{o,o})- M_N(\tilde\mu_{|n})| \le 2\,|\partial_N\! X_n|/|X_n|
\to 0 \quad(n \to \infty)\,,
$$
and the same holds for $\tilde\mu_n$. This proves (a).

To see (b), first note that $p^{(N)}(x,x) - p_{|X_n}^{(N)}(x,x) > 0$
for every $x \in \partial X_n$. The involved transition 
probabilities regard only
what happens in the ball $B(x,N/2)$, and the restriction to
$X_n$ means that only a part of that ball is admitted for the
walker, while the rest is taboo. By transitivity,
all these balls are isomorphic, and up to isometry, there are 
only finitely many ways to subdivide a ball into the admitted and taboo
parts. Hence, as $x$ varies, while $N$ is fixed, there are only finitely
many different values of $p_{|X_n}^{(N)}(x,x)$. Consequently, there
is $\ep_N > 0$ such that 
$$
p^{(N)}(x,x) - p_{|X_n}^{(N)}(x,x) \ge \ep_N \quad \forall n \in \N, 
x \in \partial X_n\,.
$$
Therefore
$$
M_N(\mu_{o,o})- M_N(\tilde\mu_{|n}) \ge \ep_N\,|\partial X_n|/|X_n|\,,
$$
which does not tend to zero as $n \to \infty$.
\end{proof}

\begin{rmk}\label{adjacency} (1) Theorem \ref{expanding}
is of course also valid  for the adjacency matrix of $X$ and its 
restriction to $X_n$, acting on $\ell^2(X)$, resp. 
$\ell^2(X_n) \subset \ell^2(X)$. \\[5pt]
(2) Part (a) can also be formulated for a sequence of finite graphs $X_n$ 
converging to $X$ that are not necessarily subgraphs of $X$. 
In that case, define 
$$
\partial^N\!X_n = \{ x \in X_n : B_{X_n}(x,N/2) \;
\mbox{is \emph{not} isomorphic with}\; B_X(o,N/2) \}\,.
$$
Then we call $(X_n)$ a \emph{F\o lner approximation} of $X$ if
$|\partial^N\!X_n|/|X_n| \to 0$ as $n \to \infty$ for every $N \in \N$.

This condition requires that  $X$ is vertex-transitive.
If $(X_n)$ is a F\o lner approximation of $X$, then the above argument 
shows that the cumulative spectral measures of $X_n$ converge  
weakly to the Plancherel measure of $X$. \\[5pt]
(3) As mentioned above, {\sc Grigorchuk and \.Zuk} 
\cite{GrZu2} consider $\Z_2 \wr \Z$, 
i.e., $DL(2,2)$, and use an approximation by Schreier graphs and the
associated cumulative spectral measures. They show that the latter converge
to the Plancherel measure. This can also be interpreted in terms of 
a F\o lner approximation. Indeed, the graphs $X_n$ defined by the action 
of the automaton considered in \cite{GrZu2} can be defined alternatively 
as follows: let $m$ be the
smallest power of $2$ that is $\ge n$, and consider the group
$\Ga_n=\Z_2^n\rtimes\Z_m$, with $\Z_m$ acting as a matrix with $1$'s on the
diagonal and just above, $0$'s elsewhere. This is the permutation group
acting on level $n$. The generator $a$ generates $\Z_m$, and $b=\sigma a$
where $\sigma$ is any non-trivial element of $\Z_2^m$. The Schreier graph
is the homogeneous space $X_n=\Z_m \backslash \Ga_n$.

The elements in $X_n$ can be naturally identified with vectors in
$\Z_2^n$. Fix an integer $N$, and let $Y_n$ denote those vectors that are
not periodic of period less than $N$. On one hand, $|Y_n|/|X_n|\to 1$
as $n\to\infty$, and on the other hand, any ball of radius $N$ around $y\in
Y_n$ embeds in the lamplighter group $\Z_2^\infty\rtimes\Z$; therefore the
graphs $X_n$ have the F\o lner approximation property.\hfill$\square$
\end{rmk}

Our computation of the spectrum of SRW on $DL(q,r)$
in \S \ref{tetrahedra} is linked with tetrahedra. In this context,
it is natural to take for $X_n$ an increasing family of tetrahedra 
$S_n$ with
height $n \to \infty\,$, whose union is $DL$, and consider the 
associated cumulative measures $\tilde \mu_{|n}$. 

If $r=q$, i.e., for the lamplighter group, we have
$|S_n| = (n+1)q^n$ and $|\partial S_n| = 2q^n$. Therefore $(S_n)$ is
a very natural F\o lner sequence in $\DL(q,q)$, which is indeed a Cayley graph
of an amenable group.

On the other hand, if $r > q$ then $|S_n| \sim r^{n+1}/(r-q)$
and $|\partial S_n| \sim r^n$, as $n \to \infty$, whence
$(S_n)$ is an expanding sequence. We know that in this case
$\DL(q,r)$ is a nonamenable graph \cite{SCWo}.

\begin{cor}\label{DLexpand}
If $r=q$ then for $S_n$, both sequences $\tilde\mu_n$ and 
$\tilde\mu_{|n}$ converge weakly to the Plancherel measure of 
$DL(q,q)$.\\[3pt]
If $r \ne q$ then $\tilde\mu_{|n}$ does not converge to the Plancherel
measure of $DL(q,r)$. 
\end{cor}

We remark that we did not compute the actual limit of 
$\tilde\mu_{|n}$, when $r \ne q$. This can be done along the lines
of the following computations.

\smallskip

Let us now consider the sequence $\tilde\mu_n$
corresponding to SRW on the graphs $S_n$ for $r \ne q$.

\emph{We shall always suppose that $r > q$.}

Let $S = S(a_1,a_2)$ be any tetrahedron with height $n$. 
We want to compute  the cumulative spectral measure of SRW $P_S$ on $S$.
We define 
$$
\partial^1\! S = \{ b_1a_2 : b_1 \in \T_q\,,\;b_1a_2 \in DL \}
\AND
\partial^2\! S = \{ a_1b_2 : b_2 \in \T_r\,,\;a_1b_2 \in DL \}
$$
(the lower and upper parts of $\partial S$). Then $P_S$ coincides
with $P$ in the interior $S \setminus \partial S$ of $S$, while
$$
\begin{aligned}
p_S(b_1a_2, x_1x_2) &= 1/r \quad\mbox{if} \quad b_1a_2 
\in \partial^1\! S\,,\; x_1 = b_1^-\,,\;x_2^-=a_2\,,\AND\\
p_S(a_1b_2, x_1x_2) &= 1/q \quad\mbox{if} \quad a_1b_2 
\in \partial^2\! S\,,\; x_1^- = a_1\,,\;x_2=b_2^-\,.
\end{aligned}
$$
In order to compare with $P$ acting on $\ell^2(X)$, it is more natural
to consider $P_S$ as a self-adjoint operator on $\ell^2(S,\mm_S)$,
where 
$$
\mm_S(x_1x_2) = \begin{cases} 
              1\,,&\mbox{if} x_1x_2 \in S \setminus \partial S\,,\\ 
              r/(r+q)\,,&\mbox{if}\; x_1x_2 \in\partial^1\! S\,,\\ 
              q/(r+q)\,,&\mbox{if}\; x_1x_2 \in\partial^2\! S\,,
              \end{cases}
$$            
instead of using the reference measure $\deg_S = (q+r)\mm_S$. 

We already know part of the spectrum of $P_S$, namely, the spectrum of
$P_S$ acting on the space of all horizontal functions with support in $S$.
We need further eigenfunctions besides the horizontal ones.
Recall the  functions $f_k^{[S,i,j]}$ and $g_m^{[S,i,j]}$
constructed in \eqref{fkgm}, $i=1, \dots, q-1$, $j=1, \dots, r-1$, with
$k,m = 1, \dots, n-1$. We can also include $k=0,n$ with 
$f_0^{[S,i,j]} = f_n^{[S,i,j]} = 0$. 
We shall now extend the range of $i$ and $j$,
adding also the values $i=0$ and $j=0$. Namely, in analogy with
\eqref{f1f2} and \eqref{fk}, and using \eqref{phii}, we define for 
$k=0, \dots, n$ and $x=x_1x_2 \in S$
\begin{equation}\label{morefk}
\begin{aligned}
f_k^{[S,0,0]}(x_1x_2) &= \begin{cases} q^{-k/2}\,r^{-(n-k)/2} 
                       &\text{if}\; x_1 \in H^1_{\hor(a_1)+k}\,,\\
                     0 & \text{otherwise\,,} 
                    \end{cases}
\\[4pt]
f_k^{[S,0,j]}(x_1x_2) &= \begin{cases} q^{-k/2}\,f_{n-k}^2[\varphi_j^2](x_2)
                       &\text{if}\; x_1 \in H^1_{\hor(a_1)+k}\,,\\
                     0 & \text{otherwise} \qquad\quad(j=1, \dots, r-1),
                    \end{cases} 
\\[4pt]
f_k^{[S,i,0]}(x_1x_2) &= \begin{cases} f_k^1[\varphi_i^1](x_1)\,r^{-(n-k)/2}
                       &\text{if}\; x_1 \in H^1_{\hor(a_1)+k}\,,\\
                     0 & \text{otherwise}\qquad\;(i=1, \dots, q-1).
                    \end{cases} 
\end{aligned}
\end{equation}
The different functions $f_k^{[S,i,j]}$ are all orthogonal, and when
$k \in \{1, \dots, n-1\}$, then they have norm one. On the other hand,
\begin{equation}\label{norm}
\| f_0^{[S,0,j]} \|^2 = \frac{q}{r+q} \AND \| f_0^{[S,i,0]} \|^2 = \frac{r}{r+q}
\quad \mbox{in}\quad \ell^2(S,\mm_S)\,.
\end{equation}
As usual, we also think of all functions of (\ref{morefk}) as being extended 
to the whole
of $DL$, with value $0$ outside of $S$. It is also convenient to set 
$f_k^{[S,i,j]} \equiv 0$ when $k < 0$ or $k > n$. 
Let $\rho = \rho(P) = \frac{2\sqrt{qr}}{q+r}$. One has 
for all pairs $(i,j)$ 
\begin{equation}\label{psaction}
\begin{aligned}
P_Sf_k^{[S,i,j]} &= \frac{\rho}{2} 
\bigl( f_{k-1}^{[S,i,j]} + f_{k+1}^{[S,i,j]} \bigr) \quad 
\begin{cases}\forall\ k \in \Z\,,&\mbox{if}\; i, j \ne 0\,,\\
      \forall\ k \in \Z \setminus \{\pm 1\}\,,&\mbox{if}\; i = 0, j \ne 0\,,\\    
      \forall\ k \in \Z \setminus \{n \pm 1\}\,,&\mbox{if}\; i \ne 0, j = 0\,,\\
      \forall\ k \in \Z \setminus \{\pm 1, n \pm 1\}\,,&\mbox{if}\; i = j = 0\,,
\end{cases}\\
P_Sf_1^{[S,0,j]} &= 
\frac{\sqrt{r}}{\sqrt{q}}f_0^{[S,0,j]} + \frac{\rho}{2}f_2^{[S,0,j]}\,, \AND\\
P_Sf_{n-1}^{[S,i,0]} &= 
\frac{\rho}{2}f_{n-2}^{[S,i,0]} + \frac{\sqrt{q}}{\sqrt{r}}f_n^{[S,i,0]}\,.
\end{aligned}
\end{equation}

\underline{Case 1.} $i, j \ne 0$. In this case, (\ref{psaction}) was 
stated in Lemma \ref{action}, and we find the
corresponding eigenvalues and eigenfunctions as in \S \ref{tetrahedra}.\\

\underline{Case 2.} $i = j = 0$. None of $f_0^{[S,0,0]}$ and $f_n^{[S,0,0]}$ 
vanish. The action of $P_S$ on the space spanned by $f_k^{[S,0,0]}$ 
($k=0,\dots, n$) is described by the $(n+1)\times(n+1)$-matrix 
$$
M^o_{n+1} = \frac{\rho}{2} 
\begin{pmatrix} 
0             & 1      &        &        &        &\\
\frac{q+r}{q} & \ddots & \ddots &        &        &\\
              & 1      & \ddots & \ddots &        &\\
              &        & \ddots & \ddots & 1      &\\
              &        &        & \ddots & \ddots & \frac{q+r}{r}\\ 
              &        &        &        &  1     & 0
\end{pmatrix}
$$
over the indices $k = 0, \dots, n$. In general, if we have an eigenvalue 
$\lambda$ of $M^o_{n+1}$ and an associated \emph{left} eigenvector
$\psi$, written as a function on $\{0, \dots, n\}$, then we obtain
a normalized eigenfunction $g$ of $P_S$ and its norm in $\ell^2(S,\mm_S)$ by 
setting
\begin{equation}\label{recipe}
g = C \sum_{k=0}^n \psi(k)\, f_k^{[S,0,0]}\,,\quad\mbox{where}\quad C^2 = 
\frac{q}{r+q}\psi(0)^2  + \sum_{k=1}^{n-1} \psi(k)^2 + \frac{r}{r+q}\psi(n)^2
\end{equation}
(using (\ref{norm})), as in \S \ref{tetrahedra}. 
Applying this recipe, we start with the following 
eigenvalues and eigenfunctions of $M^o_{n+1}\,$, resp. $P_S\,$.
\begin{equation}\label{pmone}
\begin{alignedat}{2}
\la_{0,n} &= 1\,,\quad \psi_{0,n}(k) = \sqrt{q/r}^{\,k}\,, \quad&
g_0^{[S,0,0]}(x) &= \left( 2qr \tfrac{r^n-q^n}{r^2-q^2} \right)^{-1/2}, \AND
\\
\la_{n,n} &= -1\,,\quad \psi_{n,n}(k) = \Bigl(-\sqrt{q/r}\Bigr)^{\! k}\,, 
\quad&g_n^{[S,0,0]}(x) &= 
(-1)^{\hor(x_1)}\left( 2qr\, \tfrac{r^n-q^n}{r^2-q^2} \right)^{-1/2}, 
\end{alignedat}
\end{equation}
where $x \in S$.
Next, for $m=1, \dots, n-1$, set 
$$
\beta_{m,n}^o =  \arctan \Bigl( \frac{q+r}{q-r}\tan \frac{m}{n}\pi \Bigr)\,.
$$
Some computations yield the same eigenvalue as in (\ref{eigen})
and the following corresponding eigenfunction of  $P_S\,$:
\begin{equation}\label{gmn00}
\la_{m,n} = \rho\,\cos \frac{m}{n}\pi\,,\qquad
g_m^{[S,0,0]} = \sum_{k=0}^{n} 
\frac{\sqrt{2}}{\sqrt{n}}\sin\Bigl(\frac{km}{n}\pi +\beta_{m,n}^o\Bigr) 
\, f_k^{[S,0,0]}\,.
\end{equation}
\\

\underline{Case 3.} $i =  0\,,\;j \ne 0$. Then  $f_0^{[S,0,j]}$ does not
vanish, while $f_n^{[S,0,j]} \equiv 0$. The action of $P_S$ on the space 
spanned by $f_k^{[S,0,j0]}$  ($k=0,\dots, n-1$) is described by the 
($n\times n$)-matrix 
$$
M_n' = \frac{\rho}{2} 
\begin{pmatrix} 
0             & 1      &        &        &\\
\frac{q+r}{q} & \ddots & \ddots &        &\\
              & 1      & \ddots & \ddots &\\
              &        & \ddots & \ddots & 1\\
              &        &        & 1      & 0
\end{pmatrix}
$$
over the indices $k = 0, \dots, n-1$. Computations are slightly more
involved in this case; we present the results.

\smallskip

\underline{Case 3.A} $\;n \ge \frac{r+q}{r-q}$.

\smallskip

\underline{A.1} If $n > \frac{r+q}{r-q}$ then there is precisely one 
$\al = \al_{0,n} > 0$ that solves the equation
\begin{equation}\label{largeevs}
e^{2n\al} = \frac{q-r e^{2\al}}{q e^{2\al} - r}\,.
\end{equation}
(If $n \le \frac{r+q}{r-q}$ then there is no such solution.) We get 
\begin{equation}\label{case3A1}
\begin{alignedat}{2}
\la'_{0,n} &= \rho\,\cosh \al_{0,n}\,,&\qquad&g_0^{[S,0,j]} = C'_{0,n}\sum_{k=0}^{n-1}
\sinh\bigl((n-k)\al_{0,n}\bigr) \,f_k^{[S,0,j]}\\
\la'_{n-1,n} &= -\rho\,\cosh \al_{0,n}\,,&\qquad&g_{n-1}^{[S,0,j]} = C'_{n-1,n}
\sum_{k=0}^{n-1} (-1)^k \sinh\bigl((n-k)\al_{0,n}\bigr)\,
f_k^{[S,0,j]}\,.
\end{alignedat}
\end{equation}

\smallskip

\underline{A.2} If $n = \frac{r+q}{r-q}$ then we find,
setting $\al_{0,n} = 0$ and $\al_{n-1,n} = \pi$,
\begin{equation}\label{case3A2}
\begin{alignedat}{2}
\la'_{0,n} &= \rho = \rho\,\cos \al_{0,n}\,,
&\qquad&g_0^{[S,0,j]} = C'_{0,n}\sum_{k=0}^{n-1}
\Bigl(1-\frac{k}{n}\Bigr) f_k^{[S,0,j]}\\
\la'_{n-1,n} &= -\rho = \rho\,\cos \al_{n-1,n}\,,
&\qquad&g_{n-1}^{[S,0,j]} = C'_{n-1,n}
\sum_{k=0}^{n-1} (-1)^k \Bigl(1-\frac{k}{n}\Bigr)f_k^{[S,0,j]}\,.
\end{alignedat}
\end{equation}

In both subcases A.1 and A.2, the equation
\begin{equation}\label{cot1}
\cot(n\al) = \frac{r-q}{r+q} \cot \al\,,\qquad \al \in (0\,,\,\pi)
\end{equation}
has exactly $n-2$ distinct solutions $\al_{m,m}$, $m=1, \dots, n-2$. For each 
$\al_{m,n}$, we get
\begin{equation}\label{case3A}
\la'_{m,n} = \rho\,\cos \al_{m,n}\,,\qquad
g_m^{[S,0,j]} = C'_{m,n}\sum_{k=0}^{n-1}
\sin\bigl((n-k)\al_{m,n}\bigr) \,f_k^{[S,0,j]}\,.
\end{equation}
In (\ref{case3A1}), (\ref{case3A2}) and (\ref{case3A}), the 
normalizing constants $C'_{m,n}$  are computed as $C$
in (\ref{recipe}).

\smallskip

\underline{Case 3.B} $\;n < \frac{r+q}{r-q}$.

In this case, the equation (\ref{cot1}) has exactly $n$ distinct solutions 
$\al_{m,n}$, $m=0, \dots, n-1$. Associated with each of those $\al_{m,n}$ 
there is a solution of the form (\ref{case3A}).\\

\underline{Case 4.} $i \ne  0\,,\;j = 0$. Then  $f_n^{[S,i,0]}$ does not
vanish, while $f_0^{[S,i,0]} \equiv 0$. The action of $P_S$ on the space 
spanned by $f_k^{[S,i,0]}$ ($k=1,\dots, n$) is described by the 
($n\times n$)-matrix 
$$
M_n'' = \frac{\rho}{2} 
\begin{pmatrix} 
0      & 1      &        &        &\\
1      & \ddots & \ddots &        &\\
       & \ddots & \ddots & 1      &\\
       &        & \ddots & \ddots & \frac{q+r}{r}\\ 
       &        &        &  1     & 0
\end{pmatrix}
$$
over the indices $k = 1, \dots, n$. This is analogous to Case 3.B, exchanging
$q$ with $r$ and $k,m$ with $n-k,m-k$.   The equation
\begin{equation}\label{cot2}
\cot(n\ga) = -\frac{r-q}{r+q} \cot \ga\,,\qquad \ga \in (0\,,\,\pi)
\end{equation}
has exactly $n$ distinct solutions $\ga_{m,n}$, $m=1, \dots, n$. Associated 
with each of them, we find
\begin{equation}\label{case4}
\la''_{m,n} = \rho\,\cos \ga_{m,n}\,,\qquad
g_m^{[S,i,0]} = C''_{m,n}\sum_{k=1}^{n} \sin(k\ga_{m,n}) \,f_k^{[S,i,0]}\,,
\end{equation}
with normalizing constant $C''_{m,n}$ according to (\ref{recipe}).\\

We set
$$
\begin{aligned}
\AC_S^o &= \{ g_m^{[S,0,0]} : m=0,\dots, n \}\,,\;\\
\AC_S' &= \{ g_m^{[S,0,j]} : m=0,\dots, n-1\,,\;j=1, \dots, r-1 \}\,,\;\\
\AC_S'' &= \{ g_m^{[S,i,0]} : m=1,\dots, n\,,\;i=1, \dots, q-1 \}\,.
\end{aligned}$$

\begin{lem}\label{complement}
The orthogonal complement of the subspace of horizontal functions in
$\ell^2(S)$, where $S=S(a_1,a_2)$, is spanned by 
$$
\AC_S = 
\AC_S^o \cup \bigcup \bigl\{ \AC_{S(a_1,c_2)}' : S(a_1,c_2)  \subseteq S \bigr\}
\cup \bigcup \bigl\{ \AC_{S(c_1,a_2)}'' : S(c_1,a_2)  \subseteq S \bigr\}\,,
$$
where $c_1$ and $c_2$ vary. The functions in $\AC_S$ are all orthonormal
with respect to each other.
\end{lem}

\begin{proof} Orthonormality follows from the straightforward verification
that this is true for the corresponding functions $f_m$ in the place
of the $g_m$. To prove that $\AC_S$
spans the orthogonal complement of  horizontal functions, we  
proceed as in the proof of Proposition \ref{basis} and consider
the $k$-th level $L_k$ of $S$. The space of functions
supported in $L_k$ has dimension $q^kr^{n-k}$, whence the 
codimension of the space of horizontal functions supported in $L^k$ 
is $q^k + r^{n-k} -1$. Direct counting shows that this is precisely 
the number of functions $f_m^{[S(c_1,c_2),i,j]}$ with $S(c_1,c_2) \subseteq S$,
$c_1 = a_1$ or $c_2=a_2$, and $i=0$ or $j=0$, that do not vanish on $L_k$.
\end{proof}

Thus, our long computations lead to the following result, which also shows 
that $\tilde \mu_n$ does not converge to the Plancherel measure
when $r \ne q$.
\begin{pro}\label{tildemun}
If $r > q$, then the cumulative spectral measures $\tilde \mu_n$ associated
with (non-truncated) SRW on $S_n$ converge weakly to the measure 
$\tilde \mu + \tilde \nu$, given as follows.
$$
\begin{aligned}
\tilde\mu 
= (q-1)(r-1) \sum_{N=2}^{\infty} r^{-N} \sum_{m=1}^{N-1} \de_{\la_{m,N}}
&+ (r-q)(r-1) \sum_{2 \le N \le \frac{r+q}{r-q}} r^{-N-1} 
              \sum_{m=0}^{N-1} \de_{\la_{m,N}'}\\
&+ (r-q)(r-1) \sum_{N  > \frac{r+q}{r-q}} r^{-N-1} 
              \sum_{m=1}^{N-2} \de_{\la_{m,N}'}
\end{aligned}
$$
with $\;\la_{m,N} = \rho(P) \cos(\frac{m}{M}\pi)\;$ and 
$\;\la'_{m,N} = \rho(P) \cos \al_{m,N}\;$ given by (\ref{cot1}). 
$$
\tilde \nu = (r-q)(r-1) \sum_{N  > \frac{r+q}{r-q}} r^{-N-1} 
\Bigl( \de_{\la_{0,N}'} + \de_{-\la_{0,N}'} \Bigr)\,,
$$
with $\la_{0,N}' = \rho(P) \cosh \al_{0,N}$ given by (\ref{largeevs}).

The support of the measure $\tilde \mu$ is the interval
$[-\rho(P)\,,\,\rho(P)]$. 

The support 
$\{ \pm \la_{0,N}': N  > \frac{r+q}{r-q}\}$ of $\tilde \nu$ is  
contained in $[-1\,,\,-\rho(P)] \cup [\rho(P)\,,\,1]$. The sequence
$(\la_{0,N}')$ is strictly increasing with limit $1$.
\end{pro}
 
\begin{proof}
First of all, $|S_n| = (r^{n+1} - q^{n+1})/(r-q)$. 

If we fix $n$ and consider $N \in \{2, \dots, n\}$, then $S_n = S(a_1,a_2)$
contains $(r^{n-N+1} - q^{n-N+1})/(r-q)$ different tetrahedra $S(c_1,c_2)$ 
with height $N$. With each of those, and each $i \in \{1,\dots,q-1\}$ and
$j \in \{1,\dots,r-1\}$, we associate each of the eigenvalues $\la_{m,N}$,
$m=1,\dots, N-1$. Thus, taking into account all those tetrahedra
of height $N$, we count each $\la_{m,N}$ precisely 
$(q-1)(r-1)(r^{n-N+1} - q^{n-N+1})/(r-q)$ times. If we divide
by $|S_n|$ and let $n \to \infty$, we get the first of the three parts of
$\tilde \mu$ (with an implicit use of dominated convergence).

Also, $S_n$ contains $r^{n-N}$ different tetrahedra $S(a_1,c_2)$ with height
$N$, where only $c_2$ is allowed to vary. Associated with each of them, 
and with each $j \in \{1,\dots,r-1\}$, we have each of the 
eigenvalues $\la_{m,N}'$, $m=1, \dots, N-1$. Again, taking into account all
those tetrahedra, dividing by $|S_n|$, and letting $n \to \infty$,
we obtain the second and third parts of $\tilde \mu$, plus $\tilde \nu$.
The subdivision is according to whether $|\la_{m,N}'| \le \rho(P)$ or
$> \rho(P)$, respectively.

The contributions to the spectrum of SRW on $S_n$ that come from 
$\AC_{S(a_1,a_2)}^o$ and $\AC_{S(c_1,a_2)}''$ (where $S(c_1,a_2) \subseteq
S(a_1,a_2)$ and $c_1$ varies) vanish as $n \to \infty$, because 
$q^n/r^n \to 0$.
\end{proof} 

\section{Final observations}\label{final}

{\bf A. Random walks with drift.} 
Besides SRW on $\DL(q,r)$, it may also be instructive to consider the following
variant. If the actual position is $x=x_1x_2 \in \DL(q,r)$, then we first
toss a coin, where ``head'' comes up with probability $\al \in (0\,,\,1)$.
If head comes up, then we step at random to one among the $q$ neighbours
$y=y_1y_2$ of $x$ with $y_1^- = x_1$ (i.e., downwards in Figure 2). Otherwise, 
we step at random to one among the $r$ neighbours $y$ of $x$ with $y_2^- = x_2$.
Thus, we obtain the following generalization 
$P_{\al}$ of simple random walk on $\DL(q,r)$. 
For $x=x_1x_2, y=y_1y_2 \in \DL(q,r)$
\begin{equation}\label{random-walk}
p_{\al}(x,y) = \begin{cases} \al/q 
                   & \text{if}\; y_1^- = x_1 \;\text{and}\;y_2=x_2^-\\
                                 (1-\al)/r 
                   & \text{if}\; y_1 = x_1^- \;\text{and}\;y_2^-=x_2\\
                 0 & \text{otherwise.}
                   \end{cases}
\end{equation}
In order to interpret this in terms of a lamplighter when $r=q$, 
it is best to think
of the lamps not placed at each vertex of the two-way-infinite
path $\Z$, but at the middle of each edge. Each lamp may have
$q$ different intensities or states  ($\equiv$ elements of $\Z_q$), 
the state ``off'' 
corresponding to $0 \in \Z_q$. Only finitely many lamps may be switched on.
At each step, the lamplighter tosses his $\al$-coin. If ``head'' comes up,
he moves ``down'' (from $k$ to $k+1$) and switches the lamp on the 
transversed edge to a random state. Otherwise, he moves ``up'' (to $k-1$) 
and again switches the lamp on the transversed edge to a random state.

We remark that for all values $q, r$, the random walk $P_{\al}$ on $DL(q,r)$
may be interpreted as a lamplighter walk in an extended sense.
Imagine that on each edge of $\Z$, there is a green lamp with $q$ possible
intensities (including ``off'') \emph{plus} a red lamp with $r$ possible
intensities (including ``off''). The rule is that only finitely many lamps 
may be switched on, and in addition, if the lamplighter stands at
$k$, then all lamps between $k$ and $-\infty$ have to be in a green
state, while all lamps between $k$ and $+\infty$ must be in a red
state. The lamplighter tosses his $\al$-coin. If ``head'' comes up,
he moves ``down'' (from $k$ to $k+1$) and switches the \emph{green} lamp on 
the transversed edge to a random  state, while switching off the 
red lamp on that edge. Otherwise, he moves ``up'' (to $k-1$) 
and switches the \emph{red} lamp on the transversed edge to a random state, 
while switching off the green lamp on that edge.

SRW on $\DL(q,r)$ is $P_{\al}$ with $\al = \frac{q}{q+r}$. In order to compare
$\DL(q,r)$ with $\DL(q,q)$, it may be more natural 
to consider the same $\al$ in each case.

\smallskip

For arbitrary $\al$, we have 
\begin{equation}\label{reversible}
m_{\al}(x)p_{\al}(x,y) =  m_{\al}(y)p_{\al}(y,x)\,,\quad
\mbox{where}\quad 
m_{\al}(x) = \left( \frac{\al \,r}{(1-\al)q}\right)^{\hor(x_1)},
\end{equation}
i.e., $P_{\al}$ is \emph{$m_{\al}$-reversible}. Therefore, $P_{\al}$ acts as
a self-adjoint operator on the weighted space $\ell^2(\DL,m_{\al})$
with inner product 
$$
\langle f,g\rangle_{\al} = \sum_x f(x)g(x)\,m_{\al}(x)\,.
$$
A quick computation shows that
\begin{equation}\label{conjug}
\frac{\sqrt{m_{\al}(x)}}{\sqrt{m_{\al}(y)}}\,p_{\al}(x,y) = t_{\al}\,p(x,y)
\quad \forall\ x,y \in \DL\,,\quad\mbox{where}\quad
t_{\al} = \frac{\sqrt{4\al(1-\al)}}{\rho(P)}\,.
\end{equation}
Here, $p(x,y)$ refers to simple random walk, and $\rho(P) = 2\sqrt{qr}/(q+r)$
is the spectral radius of the latter. This can be interpreted in terms of the
Hilbert space isomorphism $T_{\al}: \ell^2(\DL,m_{\al}) \to \ell^2(\DL)$
(the latter with respect to the counting measure), where
$T_{\al}f(x) = \sqrt{m _{\al}(x)}\,f(x)$. We find that
$$
T_{\al}\,P_{\al}\,f = t_{\al}\cdot P\,T_{\al}\,f \quad 
\forall\ f \in \ell^2(\DL,m_{\al})\,.
$$

\begin{cor}\label{spec_al}
The spectrum of the operator $P_{\al}$ acting on $\ell^2(\DL,m_{\al})$ is the 
the interval $\bigl[-\sqrt{4\al(1-\al)}\,,\,\sqrt{4\al(1-\al)}\bigr]$. 
It is obtained by dilating the spectrum of $P$ acting on $\ell^2(\DL)$ 
(computed in Theorem \ref{spec}) by the factor $t_{\al}$. 
The Plancherel measure associated with $P_{\al}$ is the image of the 
Plancherel measure associated with $P$ (computed in Corollary \ref{muoo})
under this dilation.
\end{cor}

\begin{cor}\label{return_al} The return probabilities of $P_{\al}$
behave asymptotically as follows.\\[3pt]
{\rm (i)} If $r > q$ then
$$
p^{(2N)}(o,o) \sim 
A_1\,\bigl(4\al(1-\al))^{N} \,
\exp\bigl(-B_1\,N^{1/3}\bigr)\, N^{-5/6}
\qquad{as}\;\; N \to \infty\,. 
$$
\\[3pt]
{\rm (ii)} If $r=q$ then 
$$
p^{(2N)}(o,o) \sim 
\bar A_1 \,\bigl(4\al(1-\al))^{N} \, \exp\bigl(-B_1\,N^{1/3}\bigr)\, N^{1/6}
\qquad{as}\;\; N \to \infty\,. 
$$
(The constants $A_1$, $\bar A_1$ and $B_1$ are as in Theorem \ref{return}). 
\end{cor}

This direct comparison of the $P_{\al}$ underlines how surprising it
is that for $r >q$ (red and green lamps with different numbers of states)
the asymptotics scale down by a factor of $N$ with respect to the case 
$r=q$ (only one type of lamps, or equivalently, red and green lamps
with the same number of states).

We remark here that {\sc Bertacchi} \cite{Ber} has proved several basic results
for general random walks on $DL(q,r)$, without being aware that they apply, 
in particular, to random walks on $\Z_q \wr \Z$. For example, when applied to
$P_{\al}$, one gets a rate of escape theorem and a central limit theorem,
as follows.

In $Z_n$ is the random vertex at time $n$ according to $P_{\al}$, with starting
point $Z_0=o$, then
$$
\frac{d(Z_n,Z_0)}{n} \to |2\al -1| \quad \mbox{almost surely.}
$$
If $\al \ne 1/2$ then
$$
\frac{d(Z_n,Z_0) - |2\al -1|n}{\sqrt{4\al(1-\al)n}} \to {\mathcal N}(0,1)
\quad \mbox{in law,}
$$
where ${\mathcal N}(0,1)$ is the standard normal distribution.
If $\al = 1/2$ then
$$
\frac{d(Z_n,Z_0)}{\sqrt{n}} \to {\mathcal M}
\quad \mbox{in law,}
$$
where ${\mathcal M}$ is the probability distribution on $\R^+$ with density
$$
\frac{d{\mathcal M}}{dt} = \frac{1}{\sqrt{2\pi}}
\bigl(e^{-t^2/4} - e^{-t^2}\bigr)\,
{\mathbf 1}_{[0\,,\,\infty)}(t)\,.
$$
Note that $2\al-1$ and $4\al(1-\al)$ are mean and variance (respectively)
of the projected random variable $\widetilde\Pi(Z_1)$ on $\Z$ (see definition
of $\widetilde\Pi$ a few lines below).
Also note that these results, contrary to Corollary \ref{return_al},
do \emph{not} differ when $r =q$, resp. $r \ne q$.\\

{\bf B. Projections.} 
We have the natural projections $\Pi_1$, $\Pi_2$ and
$\widetilde\Pi$ of $\DL(q,r)$ onto $\T_q$ and $\T_r$ respectively, namely,
$\Pi_i(x_1x_2) = x_i$ and 
$\widetilde\Pi(x_1x_2) = \hor(x_1) = \hor \circ\Pi_1(x)$.
Associated with them, we have the projected random walks 
$P_{\al,q}$, $P_{1-\al,r}$ and $\widetilde P_{\al}$, where
$$
p_{\al,q}(x_1,y_1) = \begin{cases} \al/q\,,&\mbox{if}\;y_1^- = x_1\,,\\
                                   1-\al\,,& \mbox{if}\;y_1 = x_1^-\,,\\  
                                    0   \,,& \mbox{otherwise.}
                     \end{cases}
\AND                                
\widetilde P_{\al}(k,l) = \begin{cases} \al\,,&\mbox{if}\;l=k+1\,,\\
                                   1-\al\,,& \mbox{if}\;l = k-1\,,\\  
                                    0   \,,& \mbox{otherwise.}
                     \end{cases}
$$
($P_{1-\al,r}$ is analogous to $P_{\al,q}$.) The projections are compatible
with these transition operators in the sense of factorization of Markov chains,
i.e.,
\begin{equation}\label{compatible}
\begin{aligned}
p_{\al,q}(x_1,y_1) &= \sum_{\Pi_1(y)=y_1} p_{\al}(x,y) \quad 
\mbox{for every $x$ with $\Pi_1(x)=x_1$}\,,\AND\\
\widetilde p(k,l) &= \sum_{\widetilde\Pi_1(y)=l} p_{\al}(x,y) \quad 
\mbox{for every $x$ with $\widetilde\Pi(x)=k$.}
\end{aligned}
\end{equation}
Now, $P_{\al,q}$ is reversible with respect to 
$m_{\al,q}(x_1) = \bigl(\frac{\al}{q(1-\al)}\bigr)^{\hor(x_1)}$, where 
$x_1 \in \T_q$. Also, $\widetilde P_{\al}$ is reversible with respect to 
$\tilde m_{\al}(k)= \bigl(\frac{\al}{(1-\al)}\bigr)^k$.
These operators are self-adjoint on the respective spaces 
$\ell^2(\T_q,m_{\al,q})$ and $\ell^2(\Z,\tilde m_{\al})$.
Their spectra are both known to be the interval 
$\bigl[-\sqrt{4\al(1-\al)}\,,\,\sqrt{4\al(1-\al)}\bigr]$. 
The respective Plancherel measures are continuous with respect
to Lebesgue measure on that interval; they are
\begin{equation}\label{specmeas}
\begin{aligned}
&\frac{q+1}{2\pi} \frac{\sqrt{4\al(1-\al)-\la^2}}{\tau_{\al}^2-\la^2}\,d\la
\quad \mbox{with} \;\; \tau_{\al} = \frac{\sqrt{4\al(1-\al)}}{2\sqrt{q}/(q+1)} 
\quad \mbox{for}\; P_{\al,q} \;\mbox{on}\;\T_q\,,\\[3pt]
&\frac{1}{\pi} \frac{1}{\sqrt{4\al(1-\al)-\la^2}}\,d\la
\quad \mbox{for}\; \widetilde P_{\al} \;\mbox{on}\;\Z\,.
\end{aligned}
\end{equation}
This follows from the well known results for SRW on $\T_q$ and $\Z$
(see e.g. \cite{Wbook}, p. 192 and (19.27) on p.214). To adapt the
latter methods to the $\al$-walks, one can use the same method as used above
for Corollary \ref{spec_al},
conjugating as in (\ref{conjug}) with the square root of the reversing
measure. The corresponding dilation factors are
$\tau_{\al}$ when passing to $P_{\al,q}$ from SRW on $\T_q$ (whose 
spectral radius is the denominator of $\tau_{\al}$ above), and 
$\sqrt{4\al(1-\al)}$ when passing to $\widetilde P_{\al}$ from SRW on $\Z$.

The interesting fact to observe here is that the projections preserve the
spectrum, while the respective Plancherel measures differ drastically, 
starting with a sum of point masses over a dense, countable subset, and ending
up with measures having continuous densities with respect to Lebesgue measure.\\ 

{\bf C. Green kernel estimates.} 
Using the spectral measures, one may also undertake a computation
of the asymptotic behaviour of the \emph{Green kernel}
$$
G_{\al}(o,x) = \sum_{N=0}^{\infty} p^{N}(x,y) = 
\int_{\spec(P_{\al})} \frac{1}{1-\la} \, d\mu_{o,x}(\la)\,,
\quad x, y \in DL(q,r)\,.
$$
\emph{in space}, i.e., as $d(o,x) \to \infty$. 
Due to the many oscillating terms that occur in these integrals,
the spectral method becomes quite tedious, while a more probabilistic
reasoning is more efficient and likely to admit extensions to more general
random walks, see {\sc \mbox{Brofferio} and Woess} \cite{BroWo}.
The asymptotics depend on the way (direction) in which $x$ tends to $\infty$
in $\DL$, and they also differ according to whether the drift is non-zero
($\al \ne 1/2$) or zero ($\al = 1/2$).

\end{document}